\documentclass[a4paper,10pt,reqno
]{amsart}

\usepackage[top=2.5cm, bottom=2.5cm, left=2.5cm, right=2.5cm]{geometry}

\usepackage{soul}

\usepackage{amssymb,amsmath,amsthm,ascmac,array,bm}
\usepackage[dvipdfmx]{graphicx}
\usepackage{xcolor}

\usepackage{mathtools}
\mathtoolsset{showonlyrefs} 

\newtheorem{thm}{Theorem}[section]

\newtheorem{rem}[thm]{Remark}
\newtheorem{assump}[thm]{Assumption}

\numberwithin{equation}{section}
\allowdisplaybreaks

\newcommand{\mcc}{\mathcal{C}}

\newcommand{\mci}{\mathcal{I}}

\newcommand{\mcl}{\mathcal{L}}

\newcommand{\mbbr}{\mathbb{R}}

\newcommand{\mbby}{\mathbb{Y}}

\newcommand{\mbbrp}{\mathbb{R}_{+}}

\newcommand{\al}{\alpha} \newcommand{\lam}{\lambda} \newcommand{\ep}{\epsilon} 
 \newcommand{\del}{\delta} \newcommand{\sig}{\sigma}
\newcommand{\D}{\Delta} \newcommand{\Sig}{\Sigma} \newcommand{\gam}{\gamma}
\newcommand{\Lam}{\Lambda} \newcommand{\Gam}{\Gamma}

\newcommand{\p}{\partial}
 
\newcommand{\cil}{\xrightarrow{\mcl}} 
\newcommand{\cip}{\xrightarrow{p}} 

\newcommand{\argmin}{\mathop{\rm argmin}}

\def\ds#1{\displaystyle{#1}}

\def\nn{\nonumber}

\newcommand{\pr}{P}
\newcommand{\E}{E}
\newcommand{\var}{{\rm Var}}
\newcommand{\cov}{{\rm Cov}}

\newcommand{\diag}{\mathrm{diag}}

\newcommand{\sumi}{\sum_{i=1}^{N}}

\newcommand{\tz}{\theta_{0}}
\newcommand{\tes}{\hat{\theta}_{N}}
\newcommand{\taues}{\hat{\tau}_{N}}
\newcommand{\aes}{\hat{\al}_{N}}
\newcommand{\bes}{\hat{\beta}_{N}}
\newcommand{\ges}{\hat{\gam}_{N}}
\newcommand{\les}{\hat{\lambda}_{N}}
\newcommand{\des}{\hat{\delta}_{N}}
\newcommand{\mes}{\hat{\mu}_{N}}
\newcommand{\res}{\hat{\rho}_{N}}

\def\rev#1{\textcolor{black}{#1}}

\title[Mixed-effects location-scale model based on GH distribution]
{Mixed-effects location-scale model based on generalized hyperbolic distribution}

\author[Y. Fujinaga]{Yuki Fujinaga}
\address{
Graduate School of Mathematics, Kyushu University, 744 Motooka Nishi-ku Fukuoka 819-0395, Japan
}

\author[H. Masuda]{Hiroki Masuda}
\address{Faculty of Mathematics, Kyushu University, 744 Motooka Nishi-ku Fukuoka 819-0395, Japan
\and 
Graduate School of Mathematical Sciences, The University of Tokyo, 3-8-1 Komaba Meguro-ku Tokyo 153-8914, Japan.
}
\email{hmasuda@ms.u-tokyo.ac.jp}

\date{\today}

\begin{document}
\setlength{\baselineskip}{4.5mm}

\maketitle

\begin{abstract}
Motivated by better modeling of \textit{intra-individual} variability in longitudinal data, we propose a class of location-scale mixed effects models, 
in which the data of each individual is modeled by a parameter-varying generalized hyperbolic distribution. We first study the local maximum-likelihood asymptotics and reveal the instability in the numerical optimization of the log-likelihood.
Then, we construct an asymptotically efficient estimator based on the Newton-Raphson method based on the original log-likelihood function with the initial estimator being naive least-squares-type. Numerical experiments are conducted to show that the proposed one-step estimator is not only theoretically efficient but also numerically much more stable and much less time-consuming compared with the maximum-likelihood estimator.
\end{abstract}

\section{Introduction}

The key step in the population approach \cite{Lav15} is modeling dynamics of many individuals to introduce a flexible probabilistic structure for the random vector $Y_i = (Y_i(t_{ij}))_{j=1}^{n_i} \in \mbbr^{n_i}$ representing time series data (supposed to be univariate) from $i$th individual. Here $t_{i1}<\dots<t_{in_i}$ denotes sampling times, which may vary across the individuals with possibly different $n_i$ for $i=1,\dots,N$.
The model is desired to be tractable from theoretical and computational points of view.

In the classical linear mixed-effects model \cite{LaiWar82},
the target variable $Y_i$ in $\mbbr^{n_i}$ is described by
\begin{equation}
\label{lm.model}
Y_i=X_i\beta+Z_ib_i+\epsilon_i,
\end{equation}
for $i=1,\dots,N$,
where the explanatory variables $X_i \in \mbbr^{n_i}\otimes\mbbr^{p}$ and $Z_i\in \mbbr^{n_i}\otimes\mbbr^{q}$ are known design matrices,
where $\{b_i\}$ and $\{\ep_i\}$ are mutually independent centered i.i.d. sequences with covariance matrices $G\in\mbbr^q\otimes \mbbr^q$ and $H_i\in\mbbr^{n_i}\otimes\mbbr^{n_i}$, respectively; typical examples of $H_i=(H_{i,kl})$ include $H_i=\sig^2 I_{n_i}$ ($I_q$ denotes the $q$-dimensional identity matrix) and $H_{i,kl}=\sig^2 \rho^{|k-l|}$ with $\rho$ denoting the correlation coefficient.
Although the model \eqref{lm.model} is quite popular in studying longitudinal data, it is not adequate for modeling \textit{intra-individual} variability.
Formally speaking, this means that for each $i$, conditionally on $b_i$ the objective variable $Y_i$ has the covariance which does not depend on $b_i$. Therefore the model is not suitable if one wants to incorporate a random effect across the individuals into the covariance and higher-order structures such as skewness and kurtosis.

\subsection{Mixed-effects location-scale model}
\label{sec_review.melsr}

Let us briefly review the previous study which motivated our present study.
The paper \cite{HedMerDem08} introduced a variant of \eqref{lm.model}, called the mixed-effects location-scale (MELS) model, for analyzing ecological momentary assessment (EMA) data; the MELS model was further studied in \cite{HedMerDer12}, \cite{HedDemMer09}, and \cite{HedNor13} from application and computational points of view.
\rev{
EMA is also known as the experience sampling method, which is not retrospective and the individuals are required to answer immediately after an event occurs.
Modern EMA data in mental health research is longitudinal, typically consisting of possibly irregularly spaced sampling times from each patient.
To avoid the so-called ``recall bias'' of retrospective self-reports from patients, the EMA method records many events in daily life at the moment of their occurrence.
The primary interest is modeling both between- and within-subjects heterogeneities, hence one is naturally led to incorporate random effects into both trend and scale structures.
We refer to \cite{EMA2008} for detailed information on EMA data.
}

In the MELS model, the $j$th sample $Y_{ij}$ from the $i$th individual is given by
\begin{equation}
Y_{ij} = x_{ij}^\top \beta + \exp\left(\frac12 z_{ij}^\top \al\right) \ep_{1,i} + \exp\left(\frac12 (w_{ij}^\top \tau + \sig_w \ep_{2,i})\right)\ep_{3,ij}
\label{mo.hiroki}
\end{equation}
for $1\le j\le n_i$ and $1\le i\le N$.
Here, $(x_{ij},z_{ij},w_{ij})$ are non-random explanatory variables, $(\ep_{1,i},\ep_{2,i})$ denote the i.i.d. random-effect, and $\ep_{3,ij}$ denote the driving noises for each $i\le N$ such that
\begin{equation}
(\ep_{1,i},\ep_{2,i},\ep_{3,ij}) \sim N_3\left( 0,~
\begin{pmatrix}
1 & \rho & 0 \\ \rho & 1 & 0 \\ 0 & 0 & 1
\end{pmatrix}
\right)
\nonumber
\end{equation}
and that $\ep_{3,i1}, \dots, \ep_{3,i n_i} \sim \text{i.i.d.}~N(0,1)$, with $(\ep_{1,i},\ep_{2,i})$ and $(\ep_{3,ij})_{j\le n_i}$ being mutually independent.
Direct computations give the following expressions: 
$\E[Y_{ij}]=x^{\top}_{ij}\beta$, $\var[Y_{ij}] =\exp(w^{\top}_{ij}\tau+\sigma_w^2/2) + \exp(z^{\top}_i\alpha)$, and also 
$\cov[Y_{ik},Y_{il}]=\exp(z^{\top}_i \alpha)$ for $k\ne l$; the covariance structure is to be compared with the one \eqref{GHME.cov} of our model.
Further, their conditional versions given the random-effect variable $R_i := (\ep_{1,i},\ep_{2,i})$ are as follows:
$\E[Y_{ij}|R_i] = x_{ij}^{\top}\beta+\exp(z^{\top}_i \alpha/2)\ep_{1,i}$, 
$\var[Y_{ij}| R_i] = \exp(w^{\top}_{ij}\tau + \sig_w \ep_{2,i})$, and $\cov[Y_{ik},Y_{il}| R_i] =0$ for $k\ne l$.
We also note that the conditional distribution
\begin{equation}
\mcl(Y_{i1},\dots,Y_{i n_i}| R_i) 
= N_{n_i}\left( 
X_i\beta + \mathbf{1}_{n_i}
e^{z_i^\top \al /2}
\ep_{1,i},
~\diag\big(
e^{w_{i1}^\top \tau + \sig_w \ep_{2,i}},\dots, e^{w_{i n_i}^\top \tau + \sig_w \ep_{2,i}}
\big)
\right),
\nonumber
\end{equation}
where $X_i:=(x_{i1},\dots,x_{i n_i})$ and $\mathbf{1}_{n_i}\in\mbbr^{n_i}$ has the entries all being $1$.
Importantly, the marginal distribution $\mcl(Y_{i1},\dots,Y_{i n_i})$ is not Gaussian.
See \cite{HedMerDem08} for details about the data-analysis aspects of the MELS model.

The third term on the right-hand side of \eqref{mo.hiroki} obeys a sort of normal-variance mixture with the variance mixing distribution being log-normal, introducing the so-called leptokurtosis (heavier tail than the normal distribution). Further, the last two terms on the right-hand side enable us to incorporate skewness into the marginal distribution $\mcl(Y_{ij})$; it is symmetric around $x_{ij}^\top \beta$ if $\rho=0$.

\rev{
The optimization of the corresponding likelihood function is quite time-consuming since we need to integrate the latent variables $(\ep_{1,ij},\ep_{2,ij})$:
with representing $R_i$ by the two-dimensional standard normal variable, 
the log-likelihood function of $\theta:=(\beta, \alpha, \tau, \sigma_w, \rho)$ is given by
\begin{equation}
\theta
\mapsto \sum_{i=1}^N \log \bigg\{
\int_{\mbbr^2} 
\phi_{n_i}\Big(Y_{i};\, \mu_i(\beta,\alpha,X_i, z_{i};x_1),\,
\Sig_i(\tau,\sig_w,\rho,w_{i};x_1,x_2)\Big)
\phi_{2}((x_1,x_2);0,I_2)dx_1 dx_2
\bigg\},
\label{hm:MELS-log.LF}
\end{equation}
where $w_i:=(w_{ij})_{j\le n_i}$, $z_i:=(z_{ij})_{j\le n_i}$, $\phi_m(\cdot; \mu,\Sig)$ denotes the $m$-dimensional $N(\mu,\Sig)$-density, and
\begin{align}
\mu_i(\beta,\alpha,X_i,z_{i};x_1) &:= X_i\beta + \mathbf{1}_{n_i}e^{z_i^\top \al /2} x_1,
\nn\\
\Sig_i(\tau,\sig_w,\rho,w_{i}; x_1,x_2)
&:= \diag\Big(
e^{w_{i1}^\top \tau + \sig_w (\rho x_1 + \sqrt{1-\rho^2} x_2)},\dots, e^{w_{i n_i}^\top \tau + \sig_w (\rho x_1 + \sqrt{1-\rho^2} x_2)}
\Big).
\nonumber
\end{align}
}
Just for reference, we present a numerical experiment by R Software for computing the maximum-likelihood estimator (MLE).
\rev{We set $N=1000$ and $n_1=n_2=\cdots=n_{1000}=10$} and generated $x_{ij},z_{ij},w_{ij}\sim \text{i.i.d.}~N_2(0,I_2)$ independently;
then, the target parameter is $8$-dimensional.
\rev{The true values were set as follows: $\beta=(0.6, -0.2)$, $\alpha=(-0.3,~0.5)$, $\tau=(-0.5,~0.3)$, $\sigma_w=\sqrt{0.8}\approx 0.894$, and $\rho= - 0.3$.}
The results based on a single set of data are given in Table \ref{MELSMLE}.
\rev{It took more than $20$ hours in our R code for obtaining one MLE (Apple M1 Max, memory 64GB; the R function \texttt{adaptIntegrate} was used for the numerical integration); we have also run the simulation code for $N=500$ and $n_1=n_2=\cdots=n_{500}=5$, and then it took about $8$ hours.
The program should run much faster if other software such as Fortran and MATLAB is used instead of R, but we will not deal with that direction here.
}
Though it is cheating, the numerical search started from the true values; it would be much more time-consuming and unstable if the initial values were far from the true ones.


\begin{table}[h]
\centering
\begin{tabular}{ccccccccc}
  \hline
 &$ \beta_0 $& $\beta_ 1$&$ \alpha_0 $& $\alpha_1 $&$ \tau_0 $& $\tau_1 $& $\sigma_w $& $\rho $\\ 
  \hline
True values & 0.600 & -0.200 & -0.300 & 0.500 & -0.500 & 0.300 & $0.894 $ & -0.300 \\ 
  MLE & 0.597 & -0.193 & -0.269 & 0.492 & -0.507 & 0.285 & 0.860 & -0.286 \\   
   \hline\\[-2mm]
\end{tabular}
\caption{MLE results; the computation time for one pair was about $21$ hours.}
\label{MELSMLE}
\end{table}


The EM-algorithm type approach for handling latent variables would work at least numerically, while it is also expected to be time-consuming even if a specific numerical recipe is available.
\rev{Some advanced tools for numerical integration would help to some extent, but we will not pursue it here.}

\subsection{Our objective}

In this paper, we propose an alternative computationally much simpler way of the joint modeling of the mean and within-subject variance structures.
Specifically, we construct a class of parameter-varying models based on the univariate generalized hyperbolic (GH) distribution and study its theoretical properties. The model can be seen as a special case of inhomogeneous normal-variance-mean mixtures and may serve as an alternative to the MELS model; see Section \ref{sec_gig&gh} for a summary of the GH distributions.
Recently, the family has received attention for modeling non-Gaussian continuous repeated measurement data \cite{AsaBolDigWal20}, but ours is constructed based on a different perspective directly by making some parameters of the GH distribution covariate-dependent.


This paper is organized as follows.
Section \ref{sec_model_MLE} introduces the proposed model and presents the local-likelihood analysis, followed by numerical experiments.
Section \ref{sec_optimal.estiamtor} considers the construction of a specific asymptotically optimal estimator and presents its finite-sample performance with comparisons with the MLE. Section \ref{sec_cr} gives a summary and some potential directions for related future issues.

\section{Parameter-varying generalized hyperbolic model}
\label{sec_model_MLE}


\subsection{Proposed model}

We model the objective variable at $j$th-sampling time point from the $i$th-individual by
\begin{equation}
\label{gh-melsr}
Y_{ij}=x_{ij}^{\top}\beta+s( z_{ij},\alpha)v_i+\sqrt{v_i}\, \sigma(w_{ij},\tau) \epsilon_{ij}
\end{equation}
for $j=1,\dots,n_i$ and $i=1,\dots,N$, where 
\begin{itemize}
\item $x_{ij}\in\mbbr^{p_\beta}$, $z_{ij}\in\mbbr^{p_\al'}$, and $w_{ij}\in\mbbr^{p_\tau'}$ are given non-random explanatory variables;
\item $\beta\in\Theta_\beta\subset\mathbb{R}^{p_{\beta}}$, $\alpha\in\Theta_\al\subset\mathbb{R}^{p_{\alpha}}$, and $\tau\in\Theta_\tau\subset\mathbb{R}^{p_{\tau}}$ are unknown parameters;
\item The random-effect variables $v_1,v_2,\ldots\sim\text{i.i.d.}~GIG(\lambda,\delta,\gamma)$, where GIG refers to the generalized inverse Gaussian distribution (see Section \ref{sec_gig&gh});
\item $\{\epsilon_{i}=(\epsilon_{i1},\ldots,\epsilon_{in_i})^{\top}\}_{i\ge 1}\sim\text{i.i.d.}~N(0,I_{n_i})$, independent of $\{v_i\}_{i\ge 1}$;
\item $s:\mathbb{R}^{p_{\alpha}'}\times\Theta_{\alpha}\mapsto \mathbb{R}$ and $\sigma:\mathbb{R}^{p_{\tau}'}\times\Theta_{\tau}\mapsto (0,\infty)$ are known measurable functions.
\end{itemize}
As mentioned in the introduction, for \eqref{gh-melsr} one may think of the continuous-time model without system noise:
\begin{equation}
Y_{i}(t_{ij})=x_i(t_{ij})^{\top}\beta+s( z_{i}(t_{ij}),\alpha)v_i+\sqrt{v_i}\, \sigma(w_{i}(t_{ij}),\tau) \epsilon_{i}(t_{ij}),
\nn
\end{equation}
where $t_{ij}$ denotes the $j$th sampling time for the $i$th individual.

We will write $Y_i=(Y_{i1},\ldots,Y_{in_i})\in\mbbr^{n_i}$, $x_i=(x_{i1},\ldots,x_{in_i})\in\mbbr^{n_i}\otimes\mbbr^{p_\beta}$, and so on for $i=1,\ldots,N$, and also 
\begin{equation}
\theta:=(\beta,\alpha,\tau,\lambda,\delta,\gamma) \in \Theta_\beta \times \Theta_\al \times \Theta_\tau 
\times \Theta_{\lam} \times \Theta_{\del} \times \Theta_{\gam}
=:\Theta \subset \mbbr^{p},
\nonumber
\end{equation}
where $\Theta$ is supposed to be a convex domain and $p:=p_{\beta}+p_{\alpha}+p_{\tau}+3$.
We will use the notation $(\pr_\theta)_{\theta\in\Theta}$ for the family of distributions of $\{(Y_i,v_i,\ep_i)\}_{i\ge 1}$, which is completely characterized by the finite-dimensional parameter $\theta$. The associated expectation and covariance operators will be denoted by $\E_\theta$ and $\cov_\theta$, respectively.

Let us write $s_{ij}(\al)=s(z_{ij},\alpha)$ and $\sig_{ij}(\tau)=\sigma(w_{ij},\tau)$.
For each $i\le N$, the variable $Y_{i1},\ldots,Y_{in_i}$ are $v_i$-conditionally independent and normally distributed under $\pr_\theta$:
\begin{equation}
\mcl(Y_{ij}|v_i) = N\left(x_{ij}^{\top}\beta+s_{ij}(\alpha)v_i,~\sigma^2_{ij}(\tau)v_i\right).
\nn
\end{equation}
For each $i$, we have the specific covariance structure
\begin{equation}
\cov_\theta[Y_{ij}, Y_{ik}] = \rev{s_{ij}(\alpha)s_{ik}(\alpha)}\, \var_\theta[v_i].
\label{GHME.cov}
\end{equation}
The marginal distribution $\mcl(Y_{i1},\dots,Y_{i n_i})$ is the multivariate GH distribution; a more flexible dependence structure could be incorporated by introducing the non-diagonal scale matrix (see Section \ref{sec_cr} for a formal explanation). By the definition of the GH distribution, the variables $Y_{ij}$ and $Y_{ik}$ may be uncorrelated for some $(z_{ij},\al)$ while they cannot be mutually independent.

We can explicitly write down the log-likelihood function of $(Y_1,\dots,Y_N)$ as follows:
\begin{align}
\ell_N(\theta)
&=-\frac{1}{2} \log(2\pi)\sum_{i=1}^N n_i +N\lambda\log\left(\frac{\gamma}{\delta}\right) - N\log K_\lambda(\delta\gamma) 
- \frac12 \sum_{i,j} \log\sig_{ij}^2(\tau) \nn\\
&{}\qquad 
+ 
\sum_{i=1}^N \left(\lambda-\frac{n_i}{2}\right)\log B_i(\beta,\tau,\del)
- 
\sum_{i=1}^N \left(\lambda-\frac{n_i}{2}\right)\log A_i(\al,\tau,\gam)
\nonumber\\
&{}\qquad + \sum_{i,j}
\frac{ s_{ij}(\alpha)}{\sig^2_{ij}(\tau)}(Y_{ij}-x^{\top}_{ij}\beta)
+ \sum_{i=1}^N \log K_{\lambda-\frac{n_i}{2}}\big(A_i(\al,\tau,\gam) B_i(\beta,\tau,\del)\big),
\label{log-LF}
\end{align}
where $\sum_{i,j}$ denotes a shorthand for $\sum_{i=1}^{N}\sum_{j=1}^{n_i}$ and
\begin{align}
A_i(\al,\tau,\gam) &:= \sqrt{\gamma^2 + \sum_{j=1}^{n_i}\frac{s_{ij}^2(\alpha)}{\sigma_{ij}^2(\tau)}}~,
\label{def_Ai}\\
B_i(\beta,\tau,\del) &:= \sqrt{\delta^2 + \sum_{j=1}^{n_i}\frac{1}{\sigma_{ij}^2(\tau)}(Y_{ij}-x_{ij}^{\top}\beta)^2}~.
\label{def_Bi}
\end{align}
The detailed calculation is given in Section \ref{sec_log-LF-derivation}.

To deduce the asymptotic property of the MLE, 
there are two typical ways: the global- and the local-consistency arguments.
In the present inhomogeneous model where the variables $(x_{ij},z_{ij},w_{ij})$ are non-random, the two asymptotics have different features: on the one hand, the global-consistency one generally entails rather messy descriptions of the regularity conditions as was detailed in the previous study \cite{Fuj21_m.thesis} while entailing theoretically stronger global claims; on the other hand, the local one only guarantees the existence of good local maxima of $\ell_N(\theta)$ while only requiring much weaker local-around-$\tz$ regularity conditions.

\subsection{Local asymptotics of MLE}
\label{sec_MLE.local.asymptotics}

In the sequel, we fix a true value $\tz=(\beta_0,\al_0,\tau_0,\lam_0,\del_0,\gam_0) \in \Theta$, where 
$\Theta_\del\times\Theta_\gam \subset(0,\infty)^2$; note that we are excluding the boundary (gamma and inverse-gamma) cases for $\mcl(v_i)$.

For a domain $A$, let $\mcc^k(\overline{A})$ denote a set of real-valued $\mcc^k$-class functions for which the $l$th-partial derivatives ($0\le l\le k$) admit continuous extensions to the boundary of $A$.
The asymptotic symbols will be used for $N\to\infty$ unless otherwise mentioned.

\begin{assump}
\label{MLE.as}~
\begin{enumerate}
 \item $\ds{\sup_{i\ge 1}\left( n_i \vee \max_{1\le j\le n_i}\max\{|x_{ij}|, |z_{ij}|, |w_{ij}|\} \right) < \infty}$. 
 \item $\alpha\mapsto s(z,\alpha)\in \mcc^3(\overline{\Theta_\alpha})$ for each $z$.
 \item $\tau\mapsto \sigma(w,\tau)\in \mcc^3(\overline{\Theta_\tau})$ for each $w$, and 
$\ds{\inf_{(w,\tau)\in\mathbb{R}^{p_{\tau}'}\times \Theta_{\tau}} \sigma(w,\tau)>0}$.
\end{enumerate}
\end{assump}

We are going to prove the local asymptotics of the MLE by applying the general result \cite[Theorems 1 and 2]{Swe80}.

Under Assumption \ref{MLE.as} and using the basic facts about the Bessel function $K_\cdot(\cdot)$ (see Section \ref{sec_gig&gh}), we can find a compact neighborhood $B_0\subset\Theta$ of $\tz$ such that
\begin{equation}
\forall K>0,\quad \sup_{i\ge 1}\max_{1\le j\le n_i}\sup_{\theta\in B_0} \E_\theta\big[|Y_{ij}|^K\big]<\infty.
\nonumber
\end{equation}
Note that $\min\{\del, \gam\} >0 $ inside $B_0$.

Let $M^{\otimes2} := MM^\top$ for a matrix $M$, and denote by $\lam_{\max}(M)$ and $\lam_{\min}(M)$ the largest and smallest eigenvalues of a square matrix $M$, and by $\p_\theta^k$ the $k$th-order partial-differentiation operator with respect to $\theta$. Write
\begin{equation}
\ell_N(\theta)=\sumi \zeta_i(\theta)
\nonumber
\end{equation}
for the right-hand side of \eqref{log-LF}.
Then, by the independence we have
\begin{equation}
\E_\theta\left[\left(\p_\theta\ell_N(\theta)\right)^{\otimes 2}\right] 
= \sumi \E_\theta\left[\left(\p_\theta\zeta_i(\theta)\right)^{\otimes 2}\right];
\nonumber
\end{equation}
just for reference, the specific forms of $\p_\theta\ell_N(\theta)$ and $\p_\theta^2\ell_N(\theta)$ are given in Section \ref{sec_log-LF-partial.deri}.
Further by differentiating $\theta \mapsto \p_\theta^2\ell_N(\theta)$ with recalling Assumption \ref{MLE.as}, it can be seen that
\begin{equation}
\forall K>0,\quad 
\sup_{i\ge 1}\sup_{\theta\in B_0} \E_\theta\left[ \left|\p_\theta^m\zeta_i(\theta)\right|^K \right] < \infty
\label{local.MLE_add1}
\end{equation}
for $m=1,2$, and that
\begin{equation}
\limsup_N \sup_{\theta\in B_0} \E_\theta\left[ \frac1N \sup_{\theta' \in B_0} \left|\p_\theta^3\ell_N(\theta')\right| \right] < \infty.
\label{Sweeting_check.1}
\end{equation}
These moment estimates will be used later on; unlike the global-asymptotic study \cite{Fuj21_m.thesis}, we do not need the explicit form of $\p_\theta^2\ell_N(\theta)$.

We additionally assume the diverging information condition, which is inevitable for consistent estimation:

\begin{assump}
\label{A_N.as}~
\begin{equation}
\liminf_N \inf_{\theta\in B_0} \lam_{\min}\left(
\frac{1}{N} \sumi \E_\theta\left[\left(\p_\theta\zeta_i(\theta)\right)^{\otimes 2}\right]
\right) > 0.
\nonumber
\end{equation}
\end{assump}


Under Assumption \ref{MLE.as}, we may and do suppose that the matrix
\begin{equation}
A_N(\theta) := \left(\E_\theta\left[\left(\p_\theta\ell_N(\theta)\right)^{\otimes 2}\right]\right)^{1/2}
= \left( \sumi \E_\theta\left[\left(\p_\theta\zeta_i(\theta)\right)^{\otimes 2}\right] \right)^{1/2}
\nonumber
\end{equation}
is well-defined, where $M^{1/2}$ denotes the symmetric positive-definite root of a positive definite $M$.
We also have $\sup_{\theta\in B_0}|A_N(\theta)|^{-1} \lesssim N^{-1/2}\to 0$.
This $A_N(\theta)$ will serve as the norming matrix of the MLE; see Remark \ref{rem:Studentization} below for Studentization.
Further, the standard argument through the Lebesgue dominated theorem ensures that $\E_\theta\left[\p_\theta\ell_N(\theta)\right] = 0$ and $\E_\theta\left[\left(\p_\theta\ell_N(\theta)\right)^{\otimes 2}\right]= \E_\theta\left[ - \p_\theta^2\ell_N(\theta)\right]$, followed by $A_N(\theta) = \left(\E_\theta\left[-\p^2_\theta\ell_N(\theta)\right]\right)^{1/2}$.

\medskip

For $c>0$, Assumption \ref{A_N.as} yields
\begin{align}
& \sup_{\theta':\, |\theta'-\theta|\le c/\sqrt{N}}\left|A_N(\theta)^{-1} A_N(\theta') - I_p \right|
\nn\\
&= \sup_{\theta':\, |\theta'-\theta|\le c/\sqrt{N}}
\left|\left( \frac{1}{\sqrt{N}}A_N(\theta)\right)^{-1} \left(\frac{1}{\sqrt{N}}A_N(\theta') - \frac{1}{\sqrt{N}}A_N(\theta)\right) \right|
\nn\\
&\lesssim \sup_{\theta':\, |\theta'-\theta|\le c/\sqrt{N}}
\left|\frac{1}{\sqrt{N}}A_N(\theta') - \frac{1}{\sqrt{N}}A_N(\theta)\right|
\nn\\
&\lesssim \sup_{\theta':\, |\theta'-\theta|\le c/\sqrt{N}}
\left| \left( \frac1N \sumi \E_{\theta'}\left[\left(\p_\theta\zeta_i(\theta')\right)^{\otimes 2}\right] \right)^{1/2}
- \left( \frac1N \sumi \E_\theta\left[\left(\p_\theta\zeta_i(\theta)\right)^{\otimes 2}\right] \right)^{1/2}
\right|
\nn\\
& \to 0.
\label{Sweeting_check.2}
\end{align}
Here the last convergence holds since the function $\theta\mapsto N^{-1/2}A_N(\theta)$ is uniformly continuous over $B_0$.

\medskip

Define the normalized observed information:
\begin{equation}
\mci_N(\theta) := - A_N(\theta)^{-1} \p_\theta^2\ell_N(\theta) A_N(\theta)^{-1\,\top}.
\nonumber
\end{equation}
Then, it follows from Assumption \ref{A_N.as} that
\begin{align}
\left| \mci_N(\theta)-I_p \right|
&= \left|
\left( \frac{1}{\sqrt{N}}A_N(\theta)\right)^{-1}
\left( \mci_N(\theta) - \left( \frac{1}{\sqrt{N}}A_N(\theta)\right)^{\otimes 2} \right)
\left( \frac{1}{\sqrt{N}}A_N(\theta)\right)^{-1\,\top}
\right|
\nn\\
&\lesssim \left|\mci_N(\theta) - \left( \frac{1}{\sqrt{N}}A_N(\theta)\right)^{\otimes 2} \right|
\nn\\
&\lesssim 
\left| \frac1N \sumi \left(\p_\theta^2 \zeta_i(\theta) - \E_\theta\left[\p_\theta^2 \zeta_i(\theta)\right]\right) \right|.
\nonumber
\end{align}
Then, \eqref{local.MLE_add1} ensures that
\begin{align}
\sup_{\theta\in B_0} \E_\theta\left[\left| \mci_N(\theta)-I_p \right|^2\right]
\lesssim \frac1N 
\left( \frac1N \sumi \sup_{\theta\in B_0} \E_\theta\left[\left|\p_\theta^2 \zeta_i(\theta)\right|^2\right]\right)
\lesssim \frac{1}{N}\to 0,
\nonumber
\end{align}
followed by the property
\begin{equation}
\forall\ep>0,\quad \sup_{\theta\in B_0}\pr_\theta\left[|\mci_N(\theta)-I_p|>\ep\right]\to 0.
\label{Sweeting_check.3}
\end{equation}

Let $\cil$ denote the convergence in distribution.
Having obtained \eqref{Sweeting_check.1}, \eqref{Sweeting_check.2}, and \eqref{Sweeting_check.3}, we can conclude the following theorem by applying \cite[Theorems 1 and 2]{Swe80}.


\begin{thm}
\label{hm:thm_local.MLE}
Under Assumptions \ref{MLE.as} and \ref{A_N.as}, we have the following statements under $\pr_{\tz}$.
\begin{enumerate}
\item 
For any bounded sequence $(u_N)\subset\mbbr^p$,
\begin{equation}
\ell_{N}\left(\tz+A_N(\tz)^{\top\,-1}u_N\right) - \ell_{N}\left(\tz\right) 
= u_N^{\top} \D_N(\tz) - \frac{1}{2} |u_N|^2 + o_{p}(1),
\nn
\end{equation}
with
\begin{equation}
\D_N(\tz) := A_N(\tz)^{-1} \p_{\theta}\ell_{N}(\tz) \cil N(0, I_p).
\nonumber
\end{equation}

\item 
There exists a local maximum point $\tes$ of $\ell_N(\theta)$ with $\pr_{\tz}$-probability tending to $1$, for which
\begin{equation}
A_N(\tz)^{\top}(\tes -\tz) = \D_N(\tz) + o_{p}(1) \cil N(0, I_p).
\label{hm:thm_local.MLE-1}
\end{equation}
\end{enumerate}
\end{thm}

\begin{rem}[Asymptotically efficient estimator]\normalfont
\label{rem_efficiency}
By the standard argument about the local asymptotic normality (LAN) of the family $\{\pr_\theta\}_{\theta\in\Theta}$, any estimators $\tes^\ast$ satisfying that
\begin{equation}
A_N(\tz)^{\top}(\tes^\ast -\tz) = \D_N(\tz) + o_{p}(1)
\label{def:ACe}
\end{equation}
are regular and asymptotically efficient in the sense of Haj\'{e}k-Le Cam.
See \cite{BasSco83} and \cite{Jeg82} for details.
\end{rem}

\begin{rem}[Studentization of \eqref{hm:thm_local.MLE-1}]\normalfont
\label{rem:Studentization}
\rev{Here is a remark on the construction of approximate confidence sets.}
Define the statistics
\begin{equation}
\hat{A}_N := \left(\sumi (\p_\theta\zeta_i(\tes))^{\otimes 2}\right)^{1/2}.
\label{Studentization-2}
\end{equation}
\rev{
Then, to make inferences for $\tz$ we can use the distributional approximations $\hat{A}_N(\tes-\tz) = \D_N(\tz) + o_{p}(1) \cil N_{p}(0, I_p)$ and
\begin{equation}
(\tes-\tz)^\top \hat{A}_N^2 (\tes-\tz)
\cil \chi^2(p).
\label{Studentization-3}
\end{equation}
To see this, it is enough to show that under $\pr_{\tz}$,}
\begin{equation}
A_N(\tz)^{-1}\hat{A}_N = I_p + o_p(1).
\label{Studentization-1}
\end{equation}
We have $\sqrt{N}(\tes-\tz)=O_p(1)$ by Theorem \ref{hm:thm_local.MLE} and Assumption \ref{A_N.as}.
This together with the Burkholder inequality and \eqref{local.MLE_add1} yield that
\begin{align}
N^{-1/2}\hat{A}_N &= \Bigg( 
\frac1N \sumi \left((\p_\theta\zeta_i(\tes))^{\otimes 2} - (\p_\theta\zeta_i(\tz))^{\otimes 2} \right)
\nn\\
&{}\qquad 
+ \frac1N \sumi \left((\p_\theta\zeta_i(\tz))^{\otimes 2} - \E_{\tz}\left[(\p_\theta\zeta_i(\tz))^{\otimes 2} \right]\right)
+ \left( N^{-1/2}A_N(\tz) \right)^{2} \Bigg)^{1/2}
\nn\\
&= \left( O_p(N^{-1/2}) + \left( N^{-1/2}A_N(\tz) \right)^{2} \right)^{1/2}
\nonumber
\end{align}
and hence
\begin{equation}
A_N(\tz)^{-1} \hat{A}_N = \left( N^{-1/2}A_N(\tz) \right)^{-1}
\left\{ o_p(1) + \left( N^{-1/2}A_N(\tz) \right)^{2} \right\}^{1/2}
=I_p + o_p(1),
\nonumber
\end{equation}
concluding \eqref{Studentization-1}.
Note that, instead of \eqref{Studentization-2}, we may also use 
\rev{
the square root of the observed information matrix
\begin{equation}
\widetilde{A}_N := \left(-\sumi \p_\theta^2\zeta_i(\tes)\right)^{1/2}
\label{Studentization-4}
\end{equation}
for concluding the same weak convergence as in \eqref{Studentization-3}.
In our numerical experiments, we made use of this $\widetilde{A}_N^2$ for computing the confidence interval and the empirical coverage probability.
The elements of $\widetilde{A}_N^2$ are explicit while rather lengthy: see Section \ref{sec_log-LF-partial.deri}.
}
\end{rem}

\rev{
\begin{rem}[Misspecifications]\normalfont
\label{rem_misspecification}
In addition to the linear form $x_{ij}^\top\beta$ in \eqref{gh-melsr}, misspecification of a parametric form of the function $(s(z_i,\al),\sig(w_{ij},\tau))$ is always concerned.
Using the $M$-estimation theory (for example, see \cite{Whi82} and \cite[Section 5]{Fah90}), under appropriate identifiability conditions, it is possible to handle their misspecified parametric forms.
In that case, however, the maximum-likelihood-estimation target, say $\theta_\ast$, is the optimal parameter (to be uniquely determined) in terms of the Kullback-Leibler divergence, and we do not have the LAN property in Theorem 2.3 in the usual sense while an asymptotic normality result of the form $\sqrt{N}(\tes-\theta_\ast) \cil N(0,\Gam_0^{-1}\Sig_0\Gam_0^{-1})$ could be given, where (non-random) $\Sig_0$ and $\Gam_0$ are specified by $N^{-1/2}\p_\theta\ell_N(\theta_\ast) \cil N(0,\Sig_0)$ and $-N^{-1}\p_\theta^2\ell_N(\theta_\ast) \cip \Gam_0$.
\end{rem}
}

\rev{
Finally, we note that the statistical problem will become non-standard if we allow that the true value of $(\del,\gam)$ for the GIG distribution $\mcl(v_i)$ satisfies that $\del_0=0$ or $\gam_0=0$. We have excluded these boundary cases at the beginning of Section 2.2.
}

\subsection{Numerical experiments}
\label{sec_MLE.sim}

\begin{figure}[h]
  \begin{minipage}{0.49\linewidth}
    \centering
    \includegraphics[keepaspectratio, scale=0.4]{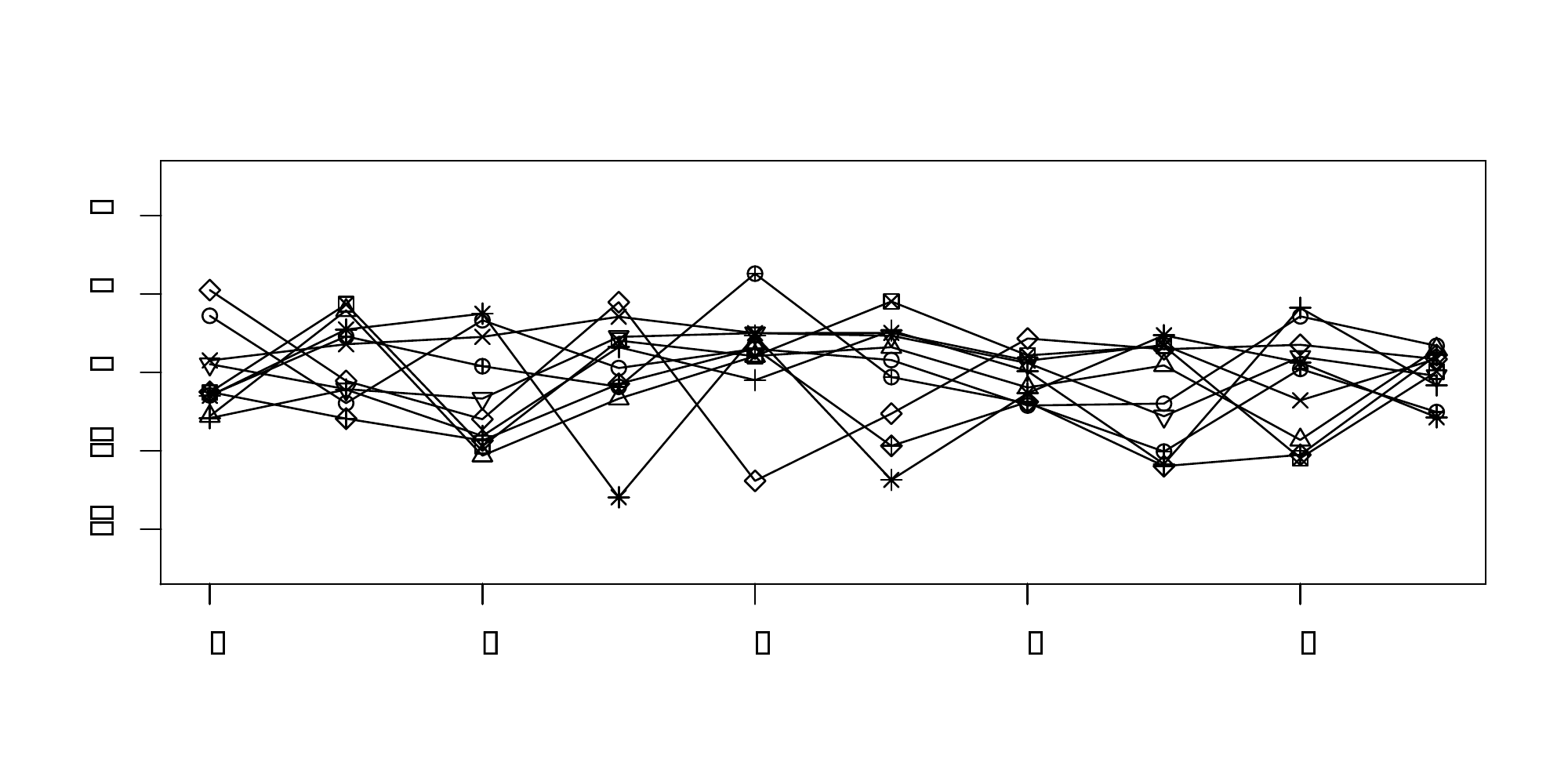}
  \end{minipage}
  \begin{minipage}{0.49\linewidth}
    \centering
    \includegraphics[keepaspectratio, scale=0.4]{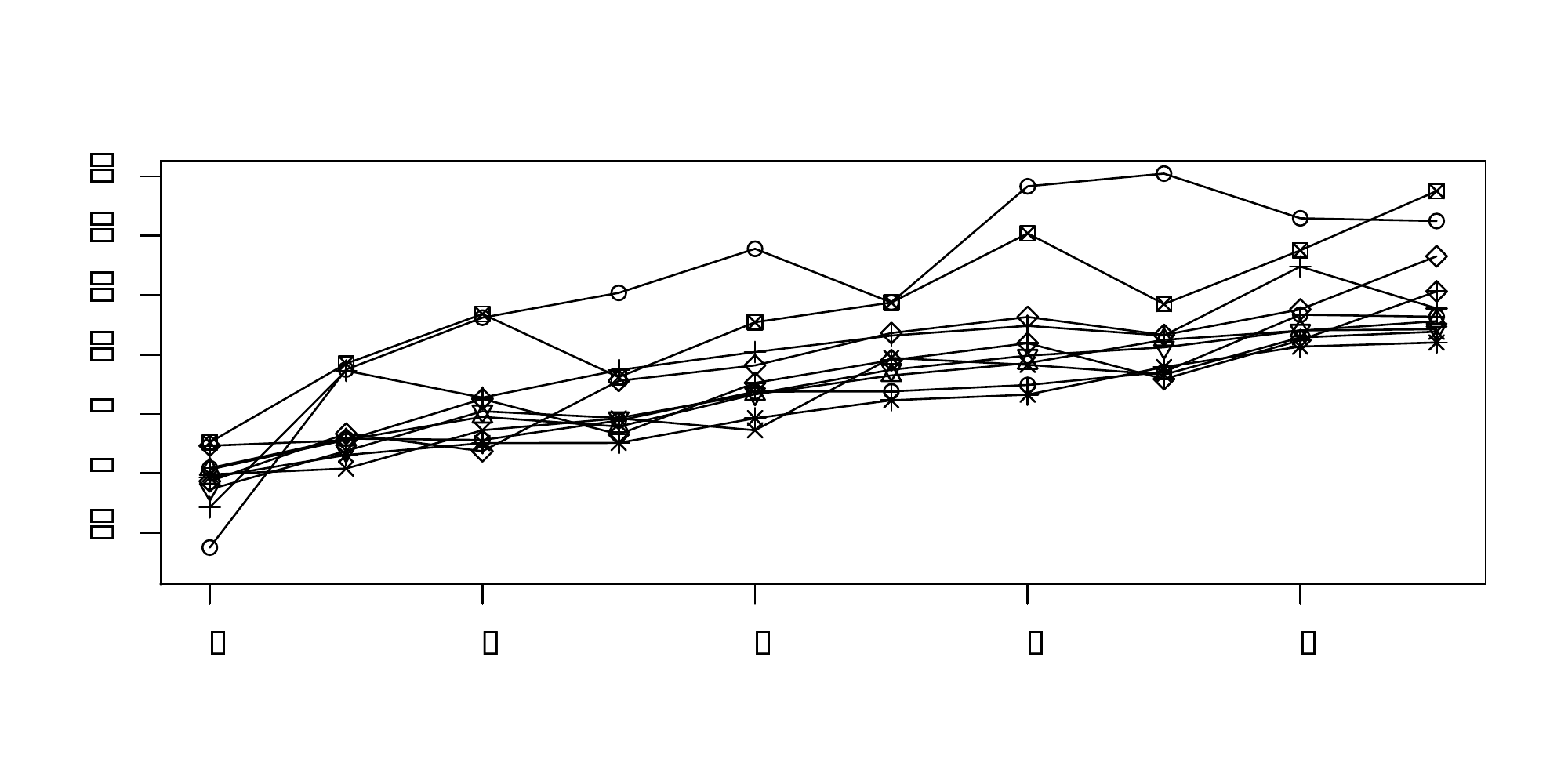}
  \end{minipage}
  \caption{Longitudinal-data plots of $10$ individuals in case (i) (left) and case (ii) (right).}
    \label{GHdata}
\end{figure}

For simulation purposes, we consider the following model:
\begin{equation}
Y_{ij}=x_{ij}^{\top}\beta+\tanh( z_{ij}^{\top}\alpha)v_i+\sqrt{v_i\exp(w_{ij}^{\top}\tau)}\,\epsilon_{ij},
\label{sim.model-1}
\end{equation}
where the ingredients are specified as follows.
\begin{itemize}
\item $N=1000$ and $n_1=n_2=\cdots =n_{1000}=10$.
\item The two different cases for the covariates $x_{ij},z_{ij},w_{ij} \in \mbbr^2$:
\begin{itemize}
\item[(i)] $x_{ij}, z_{ij}, w_{ij} \sim \text{i.i.d.}~N(0,I_2)$;
\item[(ii)] The first components of $x_{ij}, z_{ij}, w_{ij}$ are sampled from independent $N(0,1)$, and all the second ones are set to be $j-1$.
\end{itemize}
The setting (ii) incorporates similarities across the individuals; see Figure \ref{GHdata}.
\item $v_1,v_2,\dots \sim \text{i.i.d.}~GIG(\lambda,\delta,\gamma)$.
\item $\epsilon_{i}=(\epsilon_{i1},\ldots,\epsilon_{in_i})\sim N(0,I_{n_i})$, independent of $\{v_i\}$.
\item $\theta=(\beta,\alpha,\tau,\lambda,\delta,\gamma) = (\beta_0,\beta_1,\alpha_0,\alpha_1,\tau_0,\tau_1,\lambda,\delta,\gamma) \in \mbbr^9$.
\item True values of $\theta$:
\begin{itemize}
\item[(i)] $\beta=(0.3,~0.5),~\alpha=(-0.04,~0.05),~\tau=(0.05,~0.07)$, $\lambda=1.2,~\delta=1.5,~\gamma=2$;
\item[(ii)] $\beta=(0.3,~1.2),~\alpha=(-0.4,~0.8),~\tau=(0.05,~0.007)$, $\lambda=0.9,~\delta=1.2,~\gamma=0.9$.
\end{itemize}
\end{itemize}

\begin{figure}[h]
 \center
 	\includegraphics[scale=0.6]{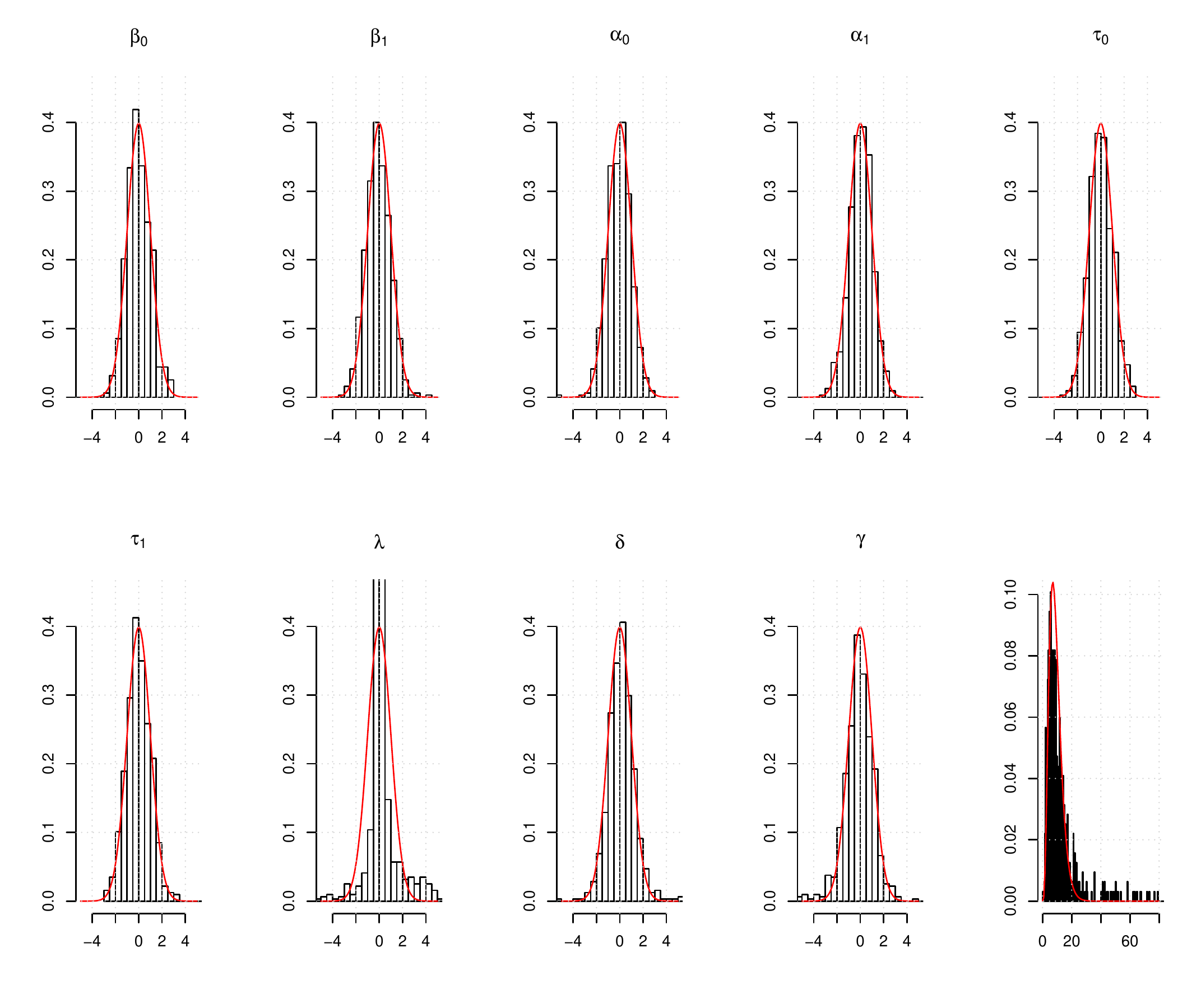}
   \caption{Standardized distributions of the MLE in case (i). The lower rightmost panel shows the chi-square approximation based on \eqref{Studentization-3}.}
      \label{MLEGH(i)}
\end{figure}

\begin{figure}[h]
 \center
  	\includegraphics[scale=0.6]{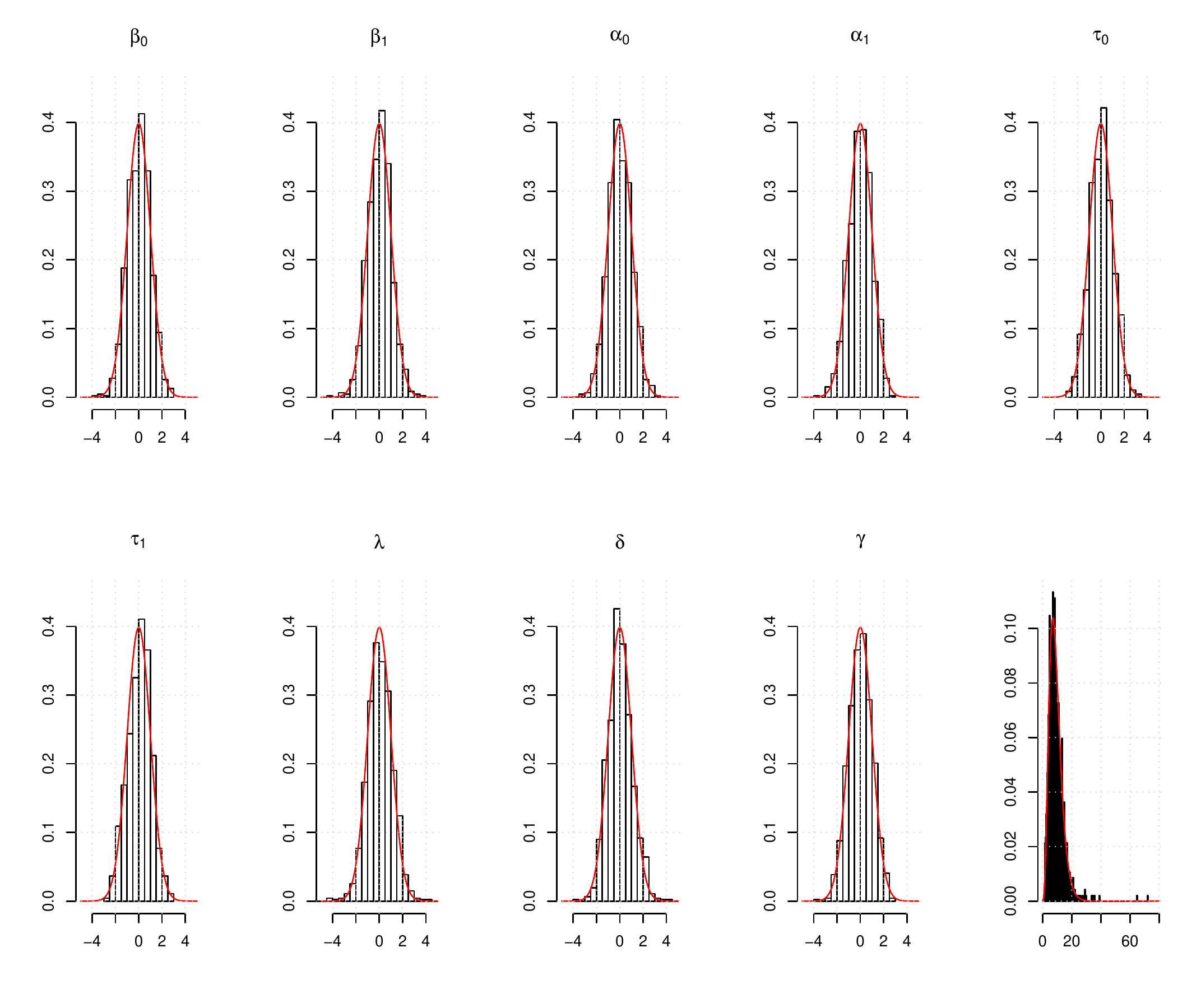}
   \caption{Standardized distributions of the MLE in case (ii). The lower rightmost panel shows the chi-square approximation based on \eqref{Studentization-3}.}
      \label{MLEGH(ii)}
\end{figure}

We numerically computed the MLE $\tes$ 
by optimizing the log-likelihood; the modified Bessel function $K_\cdot(\cdot)$ can be efficiently computed by the existing numerical libraries such as \texttt{besselK} in R Software.
We repeated the Monte Carlo trials $1000$ times, computed the 
Studentized estimates $\widetilde{A}_N(\tes-\tz)$ with \eqref{Studentization-4} in each trial, 
and then drew histograms in Figures \ref{MLEGH(i)} and \ref{MLEGH(ii)}, where the red lines correspond to the standard normal densities.
Also given in Figures \ref{MLEGH(i)} and \ref{MLEGH(ii)} are the histograms of the chi-square approximations based on \eqref{Studentization-3}.

The computation time for one MLE was about 8 minutes for case (i) and about 6 minutes for case (ii).
Estimation performance for $(\lam,\del,\gam)$ were less efficient than those for $(\beta,\al,\tau)$.
It is expected that 
the unobserved nature of the GIG variables make the standard-normal approximations relatively worse.
It is worth mentioning that case (ii) shows better normal approximations, in particular for $(\lam,\del,\gam)$; 
case (ii) would be simpler in the sense that the data from each individual have similarities in their trend (mean) structures.

\rev{
Table \ref{coverage_MLE} shows the empirical $95\%$-coverage probability for each parameter in both (i) and (ii), based on the confidence intervals $\tes^{(k)} \pm z_{\al/2}[(-\p_{\theta}^2\ell_N(\tes))^{-1}]_{kk}^{1/2}$ for $k=1,\dots,9$ with $\tes=:(\tes^{(k)})_{k\le 9}$ and $\al=0.05$.
We had $365$ and $65$ numerically unstable cases among $1000$ trials, respectively (mostly cased by a degenerate $\det(-\p_{\theta}^2\ell_N(\tes))$). Therefore, the coverage probabilities were computed based on the remaining cases.
}

\begin{table}[h]
\centering
\begin{tabular}{cccccccccc}
  \hline
 &$ \beta_0 $& $\beta_ 1$&$ \alpha_0 $& $\alpha_1 $&$ \tau_0 $& $\tau_1 $& $\lam$& $\del $& $\gam $\\ 
  \hline
Case (i) & 0.940 & 0.948 & 0.953 & 0.946 & 0.942 & 0.953 & 0.817 & 0.863 & 0.839 \\ 
Case (ii) & 0.957 & 0.952 & 0.961 & 0.954 & 0.952 & 0.949 & 0.942 & 0.945 & 0.948 \\
   \hline\\[-2mm]
\end{tabular}
\caption{The empirical $95\%$-coverage probabilities of the MLE in cases (i) and (ii) based on $1000$ trials.}
\label{coverage_MLE}
\end{table}


\medskip

Let us note the crucial problem in the above Monte Carlo trials: the objective log-likelihood is highly non-concave, hence as usual the numerical optimization suffers from the initial-value and local-maxima problems.
Here is a numerical example based on only a single set of data with $N=1000$ and $n_1=n_2=\cdots =n_{1000}=10$ as before. The same model as in \eqref{sim.model-1} together with the subsequent settings was used, except that we set $\lam=-1/2$ known from the beginning so that the latent variables $v_1,\dots,v_N$ have the inverse-Gaussian population $IG(\del,\gam)=GIG(-1/2,\del,\gam)$. For the true parameter values specified in Table \ref{GHMLE_table}, we run the following two cases for the initial values of the numerical optimization:
\begin{itemize}
\item[(i')] The true value;
\item[(ii')] $(\underbrace{1.0\times 10^{-8},\ldots,~1.0\times 10^{-8}}_{\text{$6$ times}},~1.0\times 10^{-4},~1.0\times 10^{-3})$.
\end{itemize}
\begin{table}[h]
\centering
\begin{tabular}{rrcccccccc}
  \hline
 &&$ \beta_0 $&$ \beta_1 $& $\alpha_0$&$ \alpha_1 $&$ \tau_0 $&$ \tau_1$ & $\delta$ &$ \gamma $\\ 
True value && -3.000 & 5.000 & -3.000 & 4.000 & 0.020 & -0.050 & 1.600 & 1.000 \\ 
  \hline
(i') && -3.000 & 4.999 & -3.011 & 4.017 & 0.023 & -0.047 & 1.603 & 1.005 \\ 
(ii') && -3.000 & 4.999 & -2.966 & 3.947 & 0.023 & -0.052 & 0.947 & 0.000 \\ 
   \hline
   \\[-2mm]
\end{tabular}
\caption{MLE based on single data set; the running time was about 2 minutes for case (i') and 8 minutes for case (ii');
the performance of estimating $(\del,\gam)$ in case (ii') shows instability.}
\label{GHMLE_table}
\end{table}
The results in Table \ref{GHMLE_table} clearly show that the inverse-Gaussian parameter $(\del,\gam)$ can be quite sensitive to a bad starting point for the numerical search.
In the next section, to bypass the numerical instability we will construct easier-to-compute initial estimators and their improved versions asymptotically equivalent to the MLE.

\section{Asymptotically efficient estimator}
\label{sec_optimal.estiamtor}


Building on Theorem \ref{hm:thm_local.MLE}, we now turn to global asymptotics through the classical Newton-Raphson type procedure.
A systematic account for the theory of the one-step estimator can be found in many textbooks, such as \cite[Section 5.7]{vdV98}.
Let us briefly overview the derivation with the current matrix-norming setting.

Suppose that we are given an initial estimator $\tes^0=(\aes^0,\bes^0,\taues^0,\les^0,\des^0,\ges^0)$ of $\tz$ satisfying that
\begin{equation}
\hat{u}_N^0 := A_{N}(\tz)^{\top}(\tes^0 -\tz) = O_p(1).
\nonumber
\end{equation}
By Theorem \ref{hm:thm_local.MLE} and Assumption \ref{A_N.as}, this amounts to
\begin{equation}
\sqrt{N}(\tes^0 -\tz) = O_p(1).
\label{tes0.rate}
\end{equation}
We define the one-step estimator $\tes^1$ by
\begin{equation}
\label{oseDef}
\hat{\theta}^1_N := \hat{\theta}_N^0 - \left(\partial^2_\theta\ell_N(\hat{\theta}_N^0)\right)^{-1}\partial_{\theta}\ell_N(\hat{\theta}_N^0)
\end{equation}
on the event $\{\tes^1\in\Theta,~\det(\partial^2_\theta\ell_N(\hat{\theta}_N^0))\ne 0\}$, the $\pr_{\tz}$-probability of which tends to $1$. 
Write $\hat{u}_N^1 = A_{N}(\tz)^{\top}(\tes^1 -\tz)$ and $\hat{\mci}_N^0 = -A_N(\tz)^{-1} \p^2_\theta\ell_N(\hat{\theta}_N^0) A_N(\tz)^{-1\,\top}$.
Using Taylor expansion, we have
\begin{equation}
\hat{\mci}_N^0 \hat{u}_N^1 = \hat{\mci}_N^0 \hat{u}_N^0 + A_N(\tz)^{-1} \p_\theta\ell_N(\hat{\theta}_N^0).
\label{onestep-pre1}
\end{equation}
By the arguments in Section \ref{sec_MLE.local.asymptotics}, it holds that $|\hat{\mci}_N^0| \vee |\hat{\mci}_N^{0\,-1}|=O_p(1)$.
From \eqref{tes0.rate},
\begin{align}
A_N(\tz)^{-1} \p_\theta\ell_N(\hat{\theta}_N^0) = \D_N(\tz) - \hat{\mci}_N^0 \hat{u}_N^0 + O_p\big(N^{-1/2}\big).
\label{onestep-pre2}
\end{align}
Combining \eqref{onestep-pre1} and \eqref{onestep-pre2} and recalling Remarks \ref{rem_efficiency} and \ref{rem:Studentization}, we obtain the asymptotic representation \eqref{def:ACe} for $\tes^1$, followed by the asymptotic standard normality
\begin{equation}
\hat{u}_N^1 = \D_N(\tz) + o_{p}(1) \cil N_{p}(0, I_p)
\nn
\end{equation}
and its asymptotic optimality.


\subsection{Construction of initial estimator}
\label{sec_ini.est.construction}

This section aims to construct a $\sqrt{N}$-consistent estimator $\tes^0$ satisfying \eqref{tes0.rate} through the stepwise least-squares type estimators for the first three moments of $Y_{ij}$.
We note that the model \eqref{gh-melsr} does not have a conventional location-scale structure because of the presence of $v_i$ in the two different terms.

We assume that the parameter space $\Theta_\beta \times \Theta_\al \times \Theta_\tau \times 
\Theta_\lam\times\Theta_\del\times\Theta_\gam$ is a bounded convex domain in $\mbbr^{p_\beta}\times\mbbr^{p_\al}\times\mbbr^{p_\tau}\times \mbbr\times (0,\infty)^2$ with the compact closure.
Write $\theta'=(\lam,\del,\gam)$ for the parameters contained in $\mcl(v_1)$, the true value being denoted by $\theta'_0=(\lam_0,\del_0,\gam_0)$.
Let $\mu=\mu(\theta')=\E_\theta[v_1]$, $c=c(\theta'):=\var_{\theta}[v_1]$, and $\rho=\rho(\theta'):=\E_\theta[(v_i -\E_\theta[v_i])^3]$; write $\mu_0=\mu(\theta'_0)$, $c_0=c(\theta'_0)$, and $\rho_0=\rho(\theta'_0)$ correspondingly.
Further, we introduce the sequences of the symmetric random matrices:
\begin{align}
Q_{1,N}(\al) &:= \frac1N \sum_{i,j} \big( \mu_0\,\p_\al s_{ij}(\al),\, x_{ij},\, s_{ij}(\al_0) \big)^{\otimes 2},
\nn\\
Q_{2,N}(\tau) &:= \frac1N \sum_{i,j} 
\left( \mu_0\,\p_\tau (\sig^2_{ij})(\tau),\, s_{ij}^2(\al_0) \right)^{\otimes 2}.
\nonumber
\end{align}
To state our global consistency result, we need additional assumptions.

\begin{assump}
\label{assup:initial.estimator}~
In addition to Assumption \ref{MLE.as}, the following conditions hold.
\begin{enumerate}
\item Global identifiability of $(\al,\beta,\mu)$:
\begin{enumerate}
\item $\ds{\sup_\al |Q_{1,N}(\al) - Q_1(\al) | \to 0}$ for some non-random function $Q_{1}(\al)$;
\item $\ds{\liminf_N \inf_{\al} \lam_{\min}(Q_{1,N}(\al))>0}$.
\end{enumerate}

\item Global identifiability of $(\tau,c)$:
\begin{enumerate}
\item $\ds{\sup_\tau |Q_{2,N}(\tau) - Q_2(\tau) | \to 0}$ for some non-random function $Q_{2}(\tau)$;
\item $\ds{\liminf_N \inf_{\tau}\lam_{\min}(Q_{2,N}(\tau))>0}$.
\end{enumerate}

\item Global identifiability of $\rho$: $\ds{\liminf_N \frac1N\sum_{i,j} s_{ij}^6(\al_0)>0}$.

\item 
There exists a neighborhood of $\theta'_0$ on which the mapping $\psi:\,\Theta_\lam\times\Theta_\del\times\Theta_\gam \to (0,\infty)^2\times\mbbr$ defined by $\psi(\theta')=(\mu(\theta'),c(\theta'),\rho(\theta'))$ 
is bijective, and $\psi$ is continuously differentiable at $\tz$ with nonsingular derivative.
\end{enumerate}
\end{assump}

\medskip

To construct $\tes^0$, we will proceed as follows.

\begin{description}
\item[Step 1] 
Noting that $\E_\theta[Y_{ij}]=x_{ij}^{\top}\beta+s_{ij}(\alpha) \mu$, 
we estimate $(\beta,\al,\mu)$ by minimizing
\begin{equation}
M_{1,N}(\al,\beta,\mu) := \sum_{i,j}\left( Y_{ij} - x_{ij}^{\top}\beta - s_{ij}(\alpha) \mu \right)^2.
\label{M1N_def}
\end{equation}
Let $(\aes^0,\bes^0,\mes^0)\in \argmin_{(\al,\beta,\mu)\in \overline{\Theta_\beta \times \Theta_\al \times \Theta_\mu} } M_{1,N}(\al,\beta,\mu)$.

For estimating the remaining parameters, we introduce the (heteroscedastic) residual
\begin{equation}
\hat{e}_{ij}:=Y_{ij}-x_{ij}^{\top}\bes^0 - s_{ij}(\aes^0) \mes^0,
\label{def_e.hat}
\end{equation}
which is to be regarded as an estimator of the unobserved quantity $\sqrt{v_i}\,\sig_{ij}(\tau_0)\ep_{ij}$.

\item[Step 2] 
Noting that $\var_\theta[Y_{ij}]=\sig_{ij}^2(\tau)\mu + s_{ij}^2(\al)c$, we estimate the variance-component parameter $(\tau,\al)$ by minimizing
\begin{equation}
M_{2,N}(\tau,c) := \sum_{i,j}\left( \hat{e}_{ij}^2 - \sig_{ij}^2(\tau)\mes^0 - s_{ij}^2(\aes^0)c\right)^2.
\label{M2N_def}
\end{equation}
Let $(\taues^0,\hat{c}_N^0) \in \argmin_{(\tau,c)\in \overline{\Theta_\tau} \times (0,\infty) } M_{2,N}(\tau,c)$.

\item[Step 3] 
Noting that $\E_\theta[(Y_{ij}-\E_\theta[Y_{ij}])^3] = 3 s_{ij}(\al) \sig_{ij}^2(\tau) c + s_{ij}^3(\al) \rho$, 
we estimate $\rho$ by the minimizer $\res^0$ of
\begin{equation}
M_{3,N}(\rho) := \sum_{i,j}\left( \hat{e}_{ij}^3 - 3 s_{ij}(\aes^0) \sig_{ij}^2(\taues^0) \hat{c}_N^0 - s_{ij}^3(\aes^0) \rho \right)^2,
\nonumber
\end{equation}
that is,
\begin{equation}
\res^0 := \left( \sum_{i,j} s_{ij}^6(\aes^0)\right)^{-1} \sum_{i,j} \left\{ \hat{e}_{ij}^3 - 3 s_{ij}(\aes^0) \sig_{ij}^2(\taues^0) \hat{c}_N^0\right\} s_{ij}^3(\aes^0).
\label{hat-rho}
\end{equation}

\item[Step 4]
Finally, under Assumption \ref{assup:initial.estimator}(4), we construct $\tes^{\prime 0} = (\les^0,\des^0,\ges^0)$ through the delta method by inverting $(\mes^0,\hat{c}_N^0,\res^0)$:
\begin{align}
\sqrt{N}\big(\tes^{\prime 0} - \tz'\big)
&= \sqrt{N}\left(\psi^{-1}(\mes^0,\hat{c}_N^0,\res^0) - \psi^{-1}(\mu_0,c_0,\rho_0)\right) \nn\\
&= \big(\p_{\theta'}\psi(\tz')\big)^{-1} 
\sqrt{N}\left((\mes^0,\hat{c}_N^0,\res^0) - (\mu_0,c_0,\rho_0)\right)
=O_p(1).
\nonumber
\end{align}

\end{description}

\medskip

In the rest of this section, we will go into detail about Steps 1 to 3 mentioned above and show that the estimator $\tes^0$ thus constructed satisfies \eqref{tes0.rate}; Step 4 is the standard method of moments \cite[Chapter 4]{vdV98}.

For convenience, let us introduce some notation.
The multilinear-form notation
\begin{equation}
M[u] = \sum_{i_1,\dots,i_k}M_{i_1,\dots,i_k}u_{i_1}\dots u_{i_k} \in \mbbr
\nonumber
\end{equation}
is used for $M=\{M_{i_1,\dots,i_k}\}$ and $u=\{u_{i_1},\dots u_{i_k}\}$.
For any sequence random functions $\{F_N(\theta)\}_N$ and a non-random sequence $(a_N)_N \subset (0,\infty)$, we will write $F_N(\theta)=O_p^\ast(a_n)$ and $F_N(\theta)=o_p^\ast(a_n)$ when $\sup_\theta |F_N(\theta)|=O_p(a_N)$ and $\sup_\theta |F_N(\theta)|=o_p(a_N)$ under $\pr_{\tz}$, respectively.
Further, we will denote by $m_i=(m_{i1},\dots,m_{i n_i})\in\mbbr^{n_i}$ any zero-mean (under $\pr_{\tz}$) random variables such that $m_1,\dots,m_N$ are mutually independent and $\sup_{i\ge 1}\max_{1\le j\le n_i}\E_{\tz}[|m_{ij}|^K]<\infty$ for any $K>0$; its specific form will be of no importance.

\subsubsection{Step 1}
\label{sec_1st.step}

Put $a=(\al,\beta,\mu)$ and $a_0=(\al_0,\beta_0,\mu_0)$.
By \eqref{gh-melsr} and \eqref{M1N_def}, we have
\begin{align}
\mbby_{1,N}(a)&:= \frac1N \left( M_{1,N}(a) - M_{1,N}(a_0) \right)
\nn\\
&= -\frac2N \sum_{i,j} \left( x_{ij},s_{ij}(\al_0),\mu_0\right)
\cdot \left( \beta-\beta_0, \mu-\mu_0, s_{ij}(\al) - s_{ij}(\al_0)\right) m_{ij}
\nn\\
&{}\qquad + \left(\frac1N \sum_{i,j} \left( 
\mu_0\, \p_\al s_{ij}(\tilde{\al}),\,
x_{ij},\, s_{ij}(\al_0) \right)^{\otimes 2}\right)
\left[(a-a_0)^{\otimes 2}\right]
\nn\\
&= -\frac2N \sum_{i,j} \left( x_{ij},s_{ij}(\al_0),\mu_0\right)
\cdot \left( \beta-\beta_0, \mu-\mu_0, s_{ij}(\al) - s_{ij}(\al_0)\right) m_{ij}
\nn\\
&{}\qquad + 2\left( Q_{1,N}(\tilde{\al}) - Q_{1}(\tilde{\al})\right) \left[(a-a_0)^{\otimes 2}\right]
+ 2Q_{1}(\tilde{\al}) \left[(a-a_0)^{\otimes 2}\right],
\nn
\end{align}
where $\tilde{\al}=\tilde{\al}(\al,\al_0)$ is a point lying on the segment joining $\al$ and $\al_0$.
The first term on the rightmost side equals $O_p^\ast(N^{-1/2})$.
The second term equals $o_p^\ast(1)$ by Assumption \ref{assup:initial.estimator}(1), hence we conclude that $|\mbby_{1,N}(a) - \mbby_1(a)| = o_p^\ast(1)$ for $\mbby_1(a):=2Q_{1}(\tilde{\al}) \left[(a-a_0)^{\otimes 2}\right]$.
Moreover, we have $\inf_{\al} \lam_{\min}(Q_{1}(\al))>0$ hence $\argmin \mbby_1=\{a_0\}$, followed by the consistency $\hat{a}_N \cip a_0$.

To deduce $\sqrt{N}(\hat{a}_N - a_0)=O_p(1)$ we may and do focus on the event $\{\p_a M_{1,N}(\hat{a}_N)=0\}$, on which
\begin{equation}
N^{-1}\p_a^2 M_{1,N}(\tilde{a}_N) \sqrt{N}(\hat{a}_N - a_0) = -N^{-1/2}\p_a M_{1,N}(a_0),
\label{M1.proof-2}
\end{equation}
where $\tilde{a}_N$ is a random point lying on the segment joining $\hat{a}_N$ and $\al_0$.
Observe that
\begin{equation}
-\frac{1}{\sqrt{N}}\p_a M_{1,N}(a_0) = \frac{2}{\sqrt{N}} \sum_{i,j} 
\diag\left(\p_\al s_{ij}(\al_0),\, I_{p_\beta},\, 1\right) \left[ (\mu_0, x_{ij}, s_{ij}(\al_0)) \right] m_{ij} 
=O_p(1).
\nonumber
\end{equation}
Similarly,
\begin{align}
\frac1N \p_a^2 M_{1,N}(\tilde{a}_N) &= 
-\frac{2\mu_0}{N} \sum_{i,j} \left\{ m_{ij} - 
\left( x_{ij},s_{ij}(\al_0),\mu_0\right)
\cdot \left( \tilde{\beta}_N-\beta_0, \tilde{\mu}_N-\mu_0, s_{ij}(\tilde{\al}_N) - s_{ij}(\al_0)\right) 
\right\}
\nn\\
&{}\qquad + \frac{2}{N} \sum_{i,j} 
\left(\mu_0 \p_\al s_{ij}(\tilde{\al}_N),\, x_{ij},\, s_{ij}(\al_0)\right)^{\otimes 2}.
\nonumber
\end{align}
Concerning the right-hand side, the first term equals $o_p(1)$, and the inverse of the second term does $Q_{1,N}(\tilde{\al}_N)^{-1} = \{2Q_{1,N}(\al_0) + o_p(1)\}^{-1}=O_p(1)$.
The last two displays combined with Assumption \ref{assup:initial.estimator}(1) and \eqref{M1.proof-2} conclude that $\sqrt{N}(\hat{a}_N - a_0)=O_p(1)$; it could be shown under additional conditions that $\sqrt{N}(\hat{a}_N - a_0)$ is asymptotically centered normal, while it is not necessary here.

\subsubsection{Step 2}

Write $\hat{u}_{\beta,N}=\sqrt{N}(\bes^0 -\beta_0)$, $\hat{u}_{\mu,N}=\sqrt{N}(\mes^0 -\mu_0)$, and $\hat{u}'_{\al,ij}=\sqrt{N}(s_{ij}(\aes^0) -s_{ij}(\al_0))$.
Let $b:=(\tau,c)$ and $b_0:=(\tau_0,c_0)$, and moreover
\begin{align}
e_{ij} &:= \sqrt{v_i} \sig_{ij}(\tau_0) \ep_{ij}, \nn\\
\overline{e}_{ij} &:= Y_{ij}-\E_{\tz}[Y_{ij}] =e_{ij} + s_{ij}(\al_0) (v_i-\mu_0).
\nonumber
\end{align}
We have $\hat{e}_{ij} = e_{ij} - N^{-1/2}\hat{H}_{ij}$ with $\hat{H}_{ij} := x_{ij}^\top\hat{u}_{\beta,N} + s_{ij}(\aes^0) \hat{u}_{\mu,N} + \mu_0 \hat{u}'_{\al,ij}$. 
Introduce the zero-mean random variables $\eta_{ij}:=\overline{e}_{ij}^2 - \left(\sig_{ij}^2(\tau_0)+c_0 s_{ij}^2(\al_0)\right)$.
Then, we can rewrite $M_{2,N}(b)$ of \eqref{M2N_def} as
\begin{equation}
M_{2,N}(b) 
= \sum_{i,j} \left(\overline{\eta}_{ij}(b) + \frac{1}{\sqrt{N}} \hat{B}_{ij}\right)^2,
\nonumber
\end{equation}
where
\begin{align}
\overline{\eta}_{ij}(b) &:= \eta_{ij} - \left( (\sig_{ij}^2(\tau) - \sig_{ij}^2(\tau_0)) \mes^0 + (c-c_0) s_{ij}^2(\aes^0) \right), 
\nn\\
\hat{B}_{ij} &:= -2\hat{H}_{ij} + \frac{1}{\sqrt{N}}\hat{H}_{ij}^2 -\sig_{ij}^2(\tau_0) \hat{u}_{\mu,N} - c_0 \hat{u}'_{\al,ij}.
\nonumber
\end{align}

As in Section \ref{sec_1st.step}, we observe that
\begin{align}
\mbby_{2,N}(b) 
&:= \frac1N \left( M_{2,N}(b) - M_{2,N}(b_0) \right) \nn\\
&= O_p^\ast\left(\frac{1}{\sqrt{N}}\right) + \frac1N \sum_{i,j} \left( \overline{\eta}_{ij}^2(b) - \eta_{ij}^2 \right)
\nn
\\
&= O_p^\ast\left(\frac{1}{\sqrt{N}}\right) + \frac1N \sum_{i,j} \left( (\sig_{ij}^2(\tau) - \sig_{ij}^2(\tau_0)) \mes^0 + (c-c_0)s_{ij}^2(\aes^0) \right)^2
\nn\\
&= o_p^\ast(1) + \frac1N \sum_{i,j} \left( (\sig_{ij}^2(\tau) - \sig_{ij}^2(\tau_0)) \mu_0 + (c-c_0)s_{ij}^2(\al_0) \right)^2 \nn\\
&= o_p^\ast(1) + 2Q_{2,N}(\tilde{\tau}) \left[(b-b_0)^{\otimes 2}\right]
\nonumber
\end{align}
for some point $\tilde{\tau}=\tilde{\tau}(\tau,\tau_0)$ lying on the segment joining $\tau$ and $\tau_0$.
Thus Assumption \ref{assup:initial.estimator}(2) concludes the consistency $\hat{b}_N \cip b_0$: we have $|\mbby_{2,N}(b) - \mbby_2(b)| = o_p^\ast(1)$ with $\mbby_2(b):=2Q_{2}(\tilde{\tau}) \left[(b-b_0)^{\otimes 2}\right]$ satisfying that $\inf_{\tau} \lam_{\min}(Q_{2}(\tau))>0$, hence $\argmin \mbby_2=\{b_0\}$.

The tightness $\sqrt{N}(\hat{b}_N - b_0) = O_p(1)$ can be also deduced as in Section \ref{sec_1st.step}: it suffices to note that
\begin{align}
\frac{1}{\sqrt{N}}\p_b M_{2,N}(b_0) 
&= \frac{2}{\sqrt{N}} \sum_{i,j} \left( \eta_{ij} + \frac{1}{\sqrt{N}}\hat{B}_{ij}\right) \p_b \overline{\eta}_{ij}^2(b_0) \nn\\
&= -\frac{2}{\sqrt{N}} \sum_{i,j} \left( \mu_0\,\p_\tau (\sig^2_{ij})(\tau),\, s_{ij}(\al_0) \right) \eta_{ij} + O_p(1) = O_p(1),
\nonumber
\end{align}
and that
\begin{align}
\frac{1}{N}\p_b^2 M_{2,N}(\tilde{b}_N) = o_p(1) + 2Q_{2,N}(\tau_0)
\nonumber
\end{align}
for every random sequence $(\tilde{b}_N)$ such that $\tilde{b}_N \cip b_0$.

\subsubsection{Step 3}

By the explicit expression \eqref{hat-rho} and the $\sqrt{N}$-consistency of $(\aes^0,\bes^0,\mes^0,\hat{c}_N^0)$, we obtain
\begin{align}
& \Bigg( \frac1N \sum_{i,j} s_{ij}^6(\aes^0)\Bigg) \sqrt{N}(\res^0 - \rho_0) \nn\\
&= 
\frac{1}{\sqrt{N}} \sum_{i,j} s_{ij}^3(\aes^0)\left( \hat{e}_{ij}^3 - 3 s_{ij}(\aes^0) \sig_{ij}^2(\taues^0) \hat{c}_N^0 - s_{ij}^3(\aes^0) \rho_0 \right)
\nn\\
&= O_p(1) + 
\frac{1}{\sqrt{N}} \sum_{i,j} s_{ij}^3(\aes^0)\left\{
\left(\overline{e}_{ij} - \frac{1}{\sqrt{N}}\hat{H}_{ij}\right)^3 - 3 s_{ij}(\al_0) \sig_{ij}^2(\tau_0)c_0 - s_{ij}^3(\al_0) \rho_0 
\right\}
\nn\\
&= O_p(1) + 
\frac{1}{\sqrt{N}} \sum_{i,j} s_{ij}^3(\al_0)\left(
\overline{e}_{ij}^3 - 3 s_{ij}(\al_0) \sig_{ij}^2(\tau_0)c_0 - s_{ij}^3(\al_0) \rho_0 
\right) = O_p(1).
\nonumber
\end{align}
Hence $\sqrt{N}(\res^0-\rho_0)=O_p(1)$ under Assumption \ref{assup:initial.estimator}(3).


\medskip

We end this section with a few remarks.

\begin{rem}\normalfont
As an alternative to \eqref{M1N_def}, one could also use the profile least-squares estimator \cite{Ric61}: first we construct the explicit least-squares estimator of $(\beta,\mu)$ knowing $\al$, and then optimize $\al\mapsto M_{1,N}(\al,\bes(\al),\mes(\al))$ to get an estimator of $\al$.
\end{rem}

\begin{rem}\normalfont
If one component of $\theta'=(\lam,\del,\gam)$ is known from the very beginning, then it is enough to look at the estimation of $(\mu, c)$ and we can remove Assumption \ref{assup:initial.estimator}(3) with modifying Assumption \ref{assup:initial.estimator}(4).
\end{rem}

\begin{rem}\normalfont
Because of the asymptotic nature, the same flow of estimation procedures (the MLE, the initial estimator, and the one-step estimator) remain valid even if we replace the trend term $x_{ij}^\top \beta$ in \eqref{mo.hiroki} by some nonlinear one, say $\mu(x_{ij},\beta)$, with associated identifiability conditions.
\end{rem}

\rev{
\begin{rem}\normalfont 
We can construct a one-step estimator for the MELS model \eqref{mo.hiroki} in a similar manner to Steps 1 to 3 described in Section \ref{sec_ini.est.construction}. 
To construct an initial estimator $\tes^0=(\bes^0, \aes^0, \taues^0, \hat{\sig}^{2,0}_w, \res^0)$, we use the identities $\E_\theta[Y_{ij}]=x_{ij}^{\top}\beta$, $\var_\theta[Y_{ij}]=\exp(w^{\top}_{ij}\tau+\sigma_w^2/2) + \exp(z^{\top}_i\alpha)$, and $\E_\theta[(Y_{ij} - \E_\theta[Y_{ij}])^3]=3\sig_w \exp(z_{ij}^{\top}\al /2 + \sigma_w^2/2)\rho$.
Then, we can obtain $\bes^0$ in Step 1, $(\aes^0,\taues^0, \hat{\sig}^{2,0}_{w,N})$ in Step 2, and then $\res^0$ in Step 3 in this order through the contrast functions to be minimized: denoting $\hat{e}'_{ij} := Y_{ij}-x_{ij}^{\top}\bes^0$, we have
\begin{align}
\beta &\mapsto \sum_{i,j}\left( Y_{ij} - x_{ij}^{\top}\beta \right)^2,
\nn\\
(\al,\tau,\sig_w^2) &\mapsto \sum_{i,j}\left( \hat{e}_{ij}^{\prime \,2} - \exp(w^{\top}_{ij}\tau+\sigma_w^2/2) - \exp(z^{\top}_i\alpha) \right)^2,
\nn\\
\rho &\mapsto \sum_{i,j}\left( \hat{e}_{ij}^{\prime \, 3} - 3 \sqrt{\hat{\sig}_{w,N}^{2,0}} \, \exp(z_{ij}^{\top}\aes^0 /2 + \hat{\sig}_{w,N}^{2,0} /2)\rho
\right)^2.
\nonumber
\end{align}
As in the case of \eqref{hat-rho}, $\res^0$ is explicitly given while the meaning of the parameter $\rho$ is different in the present context.
It is also possible to develop an asymptotic theory for the MLE of the MELS and the relate one-step estimator in similar ways to the present study.
However, the one-step estimator toward the log-likelihood function \eqref{hm:MELS-log.LF} still necessitates the numerical integration over $\mbbr^2$ with respect to the two-dimensional standard normal random variables; 
the numerical integration would need to be performed for every $i=1,\dots,N$ and $j=1,\dots, n_i$, hence the computational load would still be significant.
\end{rem}
}

\subsection{Numerical experiments}


Let us observe the finite-sample performance of the initial estimator $\tes^0$, the one-step estimator $\tes^1$, and the MLE $\tes$.
The setting is as follows:
\begin{equation}
Y_{ij}=x_{ij}^{\top}\beta+\tanh( z_{ij}^{\top}\alpha)v_i+\sqrt{v_i\exp(w_{ij}^{\top}\tau)}\,\epsilon_{ij},
\label{model_ex.ose}
\end{equation}
where
\begin{itemize}
\item $N=1000$, \quad $n_1=n_2=\cdots =n_{N}=10$.
\item $x_{ij},~z_{ij},~w_{ij} \in \mbbr^2 \sim \text{i.i.d.}~N_2(0,I_2)$.
\item $v_1,v_2,\ldots\sim \text{i.i.d.}~IG(\delta,\gamma)=GIG(-1/2,\delta,\gamma)$, the inverse-Gaussian random-effect distribution.
\item $\epsilon_{i}=(\epsilon_{i1},\ldots,\epsilon_{in_i}) \sim \text{i.i.d.}~N(0,I_{n_i})$, independent of $\{v_i\}$.
\item $\theta=(\beta,\alpha,\tau,\delta,\gamma) = (\beta_0,\beta_1,\alpha_0,\alpha_1,\tau_0,\tau_1,\delta,\gamma) \in \mbbr^8$.
\item True values are $\beta=(3,5),~\alpha=(-4,5)$, $\tau=(0.05, 0.07),~\delta=1.5,~\gamma=0.7$.
\end{itemize}
In this case $\theta'=(\del,\gam)\in(0,\infty)^2$ and we need only $(\mes^0,\hat{c}_N^0)$: we have $\mu=\E_{\theta'}[v_i]=\delta/\gamma$ and $c=\var_{\theta'}[v_i]=\delta/\gamma^3$, namely
\begin{equation}
\gamma=\sqrt{\frac{\mu}{c}},\qquad \delta=\mu\gamma=\sqrt{\frac{\mu^3}{c}}.
\nonumber
\end{equation}
As initial values for numerical optimization, we set the following two different cases:
\begin{itemize}
\item[(i')] The true value;
\item[(ii')] $(1.0\times 10^{-8},\ldots,~1.0\times 10^{-8},~1.0\times 10^{-4},~1.0\times 10^{-3})$.
\end{itemize}

\medskip

\begin{figure}[h]
 \centering
   \begin{minipage}{0.49\hsize}
  \begin{center}
   \includegraphics[scale=0.38]{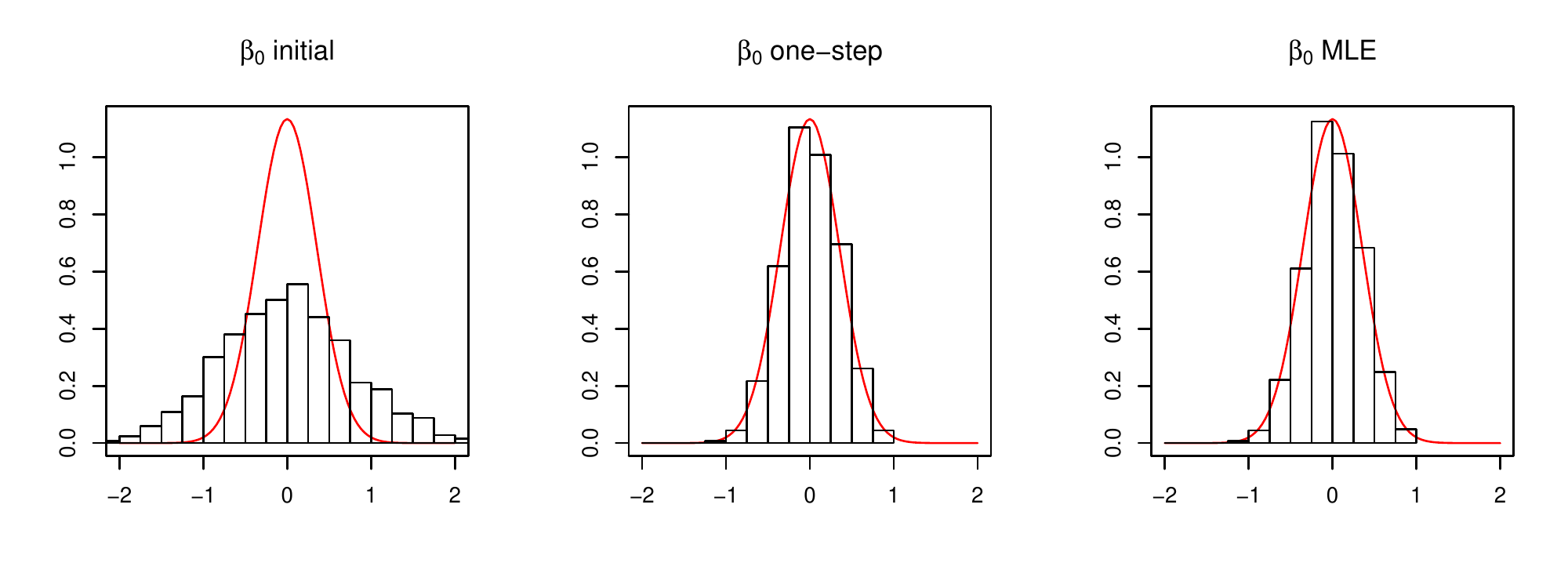}
  \end{center}
 \end{minipage}
 \begin{minipage}{0.49\hsize}
  \begin{center}
   \includegraphics[scale=0.38]{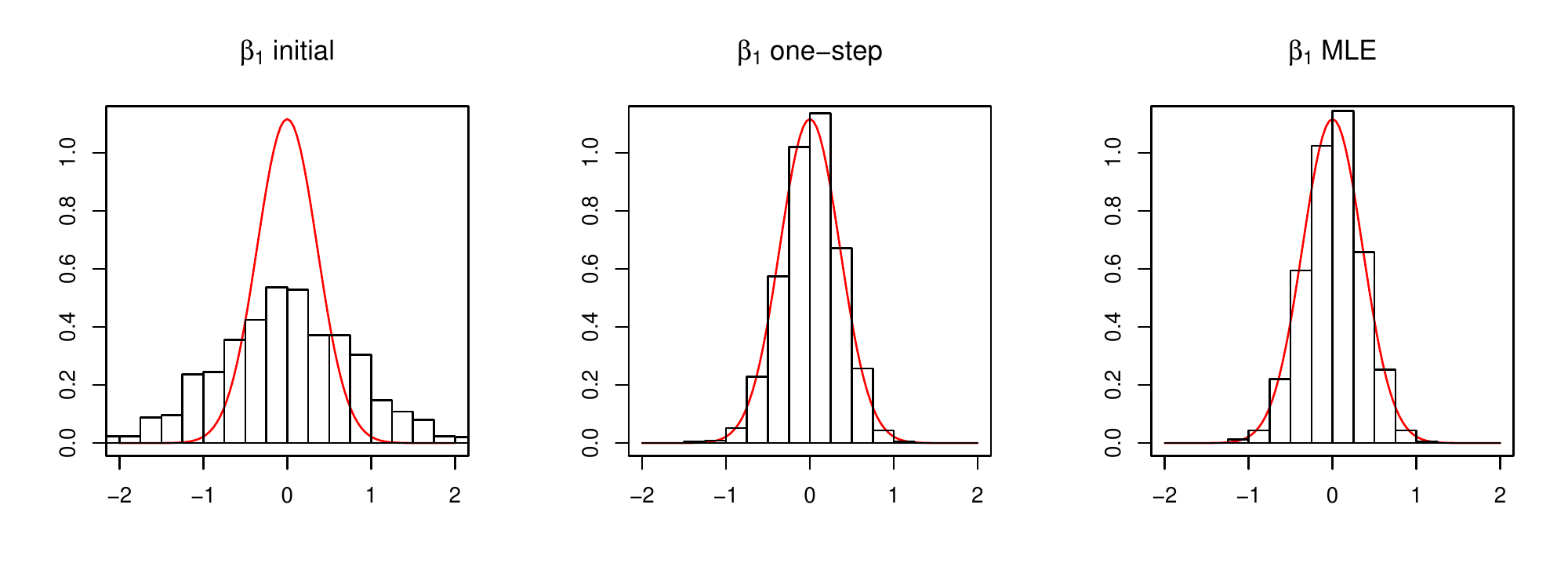}
  \end{center}
 \end{minipage}
 \\
  \begin{minipage}{0.49\hsize}
  \begin{center}
   \includegraphics[scale=0.38]{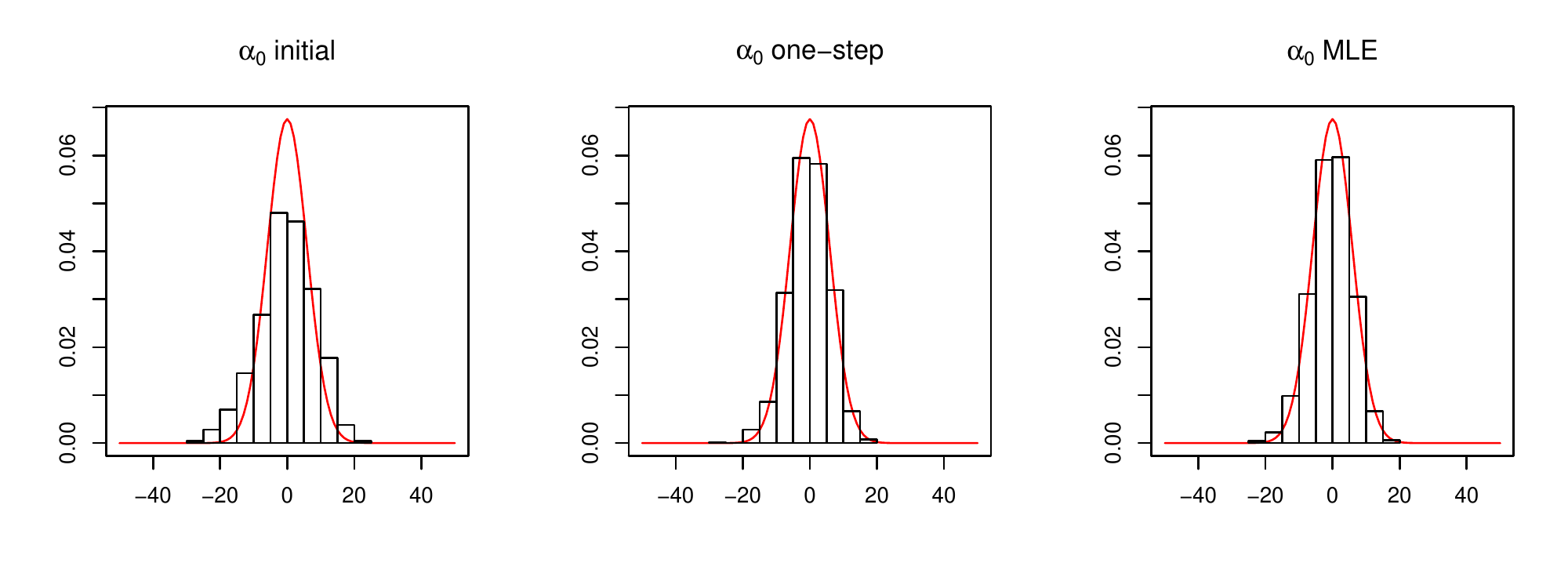}
  \end{center}
 \end{minipage}
 \begin{minipage}{0.49\hsize}
  \begin{center}
   \includegraphics[scale=0.38]{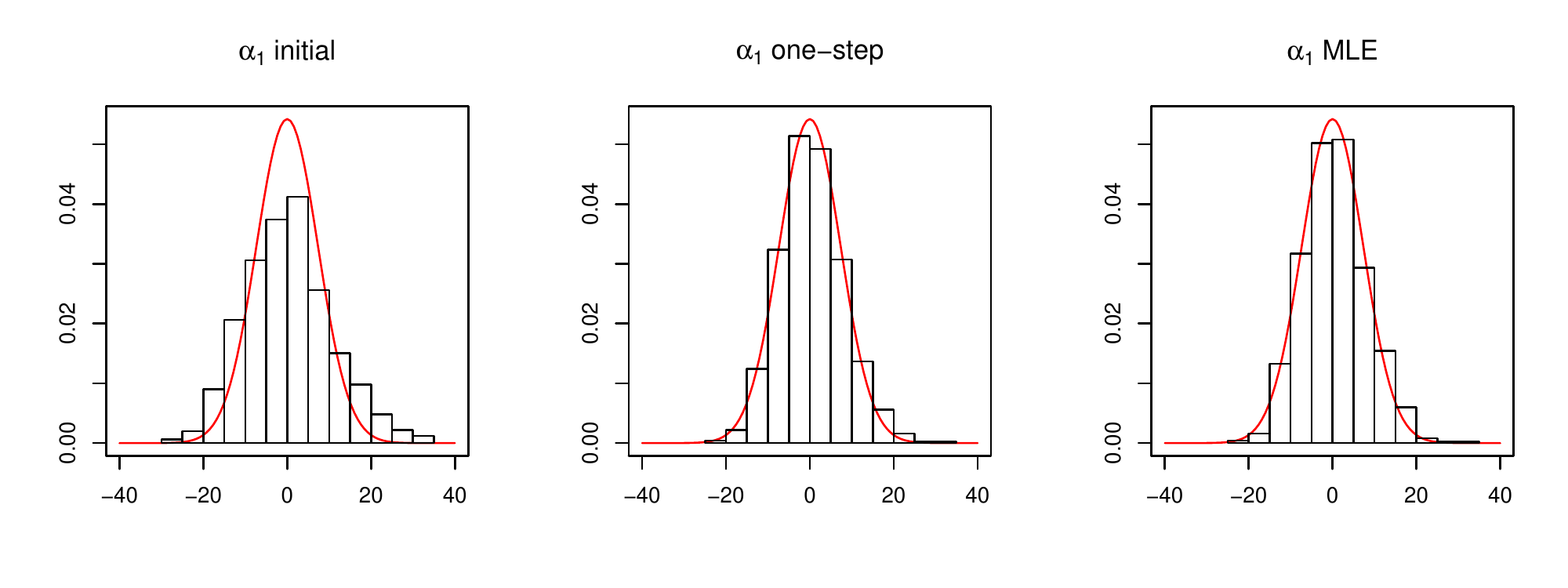}
  \end{center}
 \end{minipage}
 \\
 \begin{minipage}{0.49\hsize}
  \begin{center}
   \includegraphics[scale=0.38]{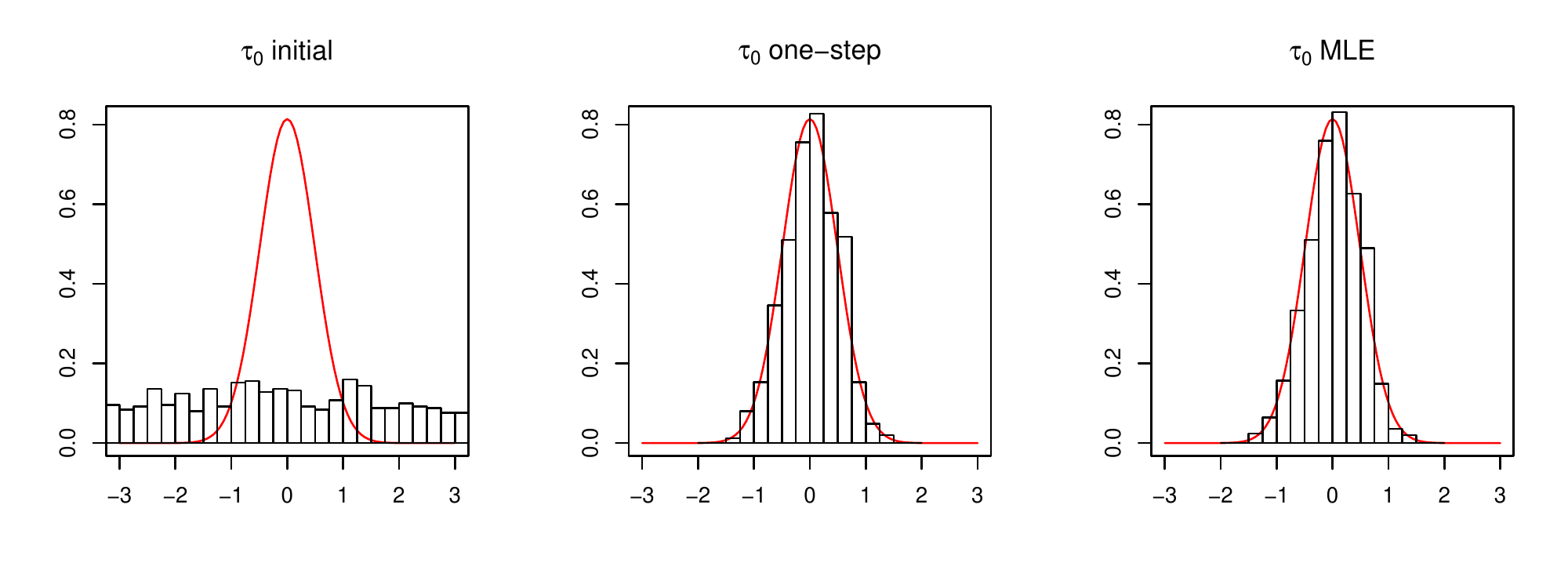}
  \end{center}
 \end{minipage}
 \begin{minipage}{0.49\hsize}
  \begin{center}
   \includegraphics[scale=0.38]{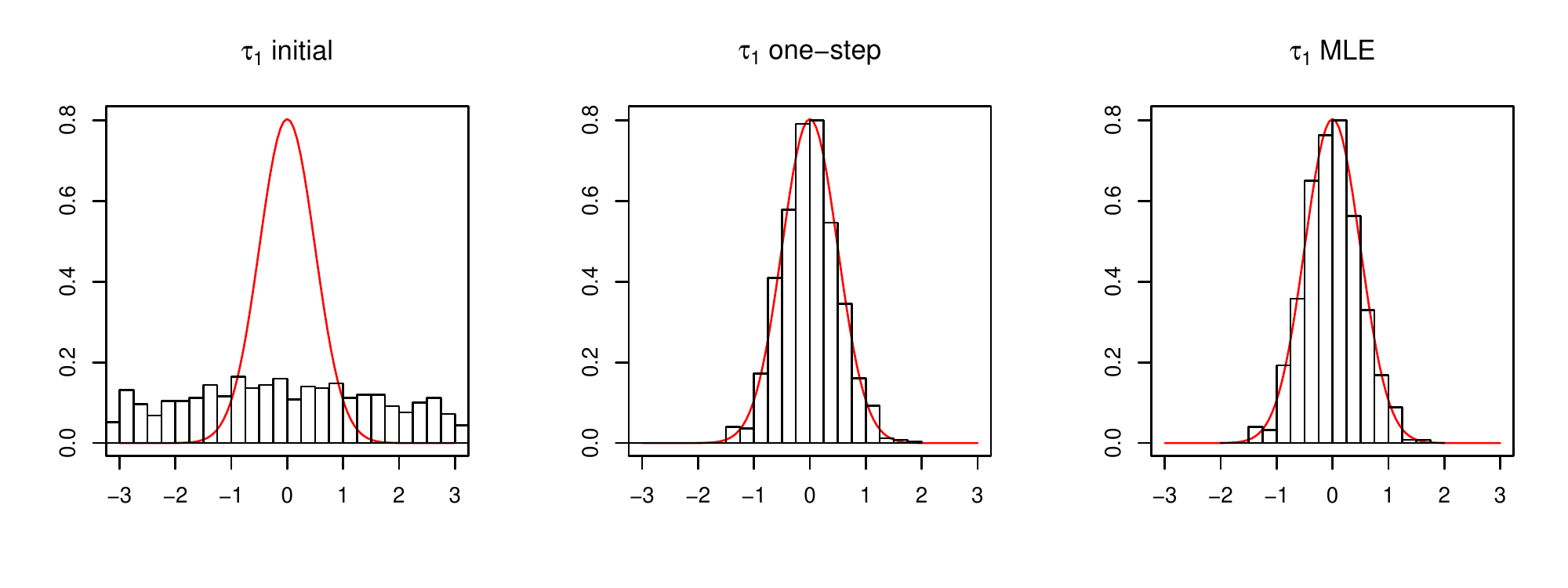}
  \end{center}
 \end{minipage}
 \\
  \begin{minipage}{0.49\hsize}
  \begin{center}
   \includegraphics[scale=0.38]{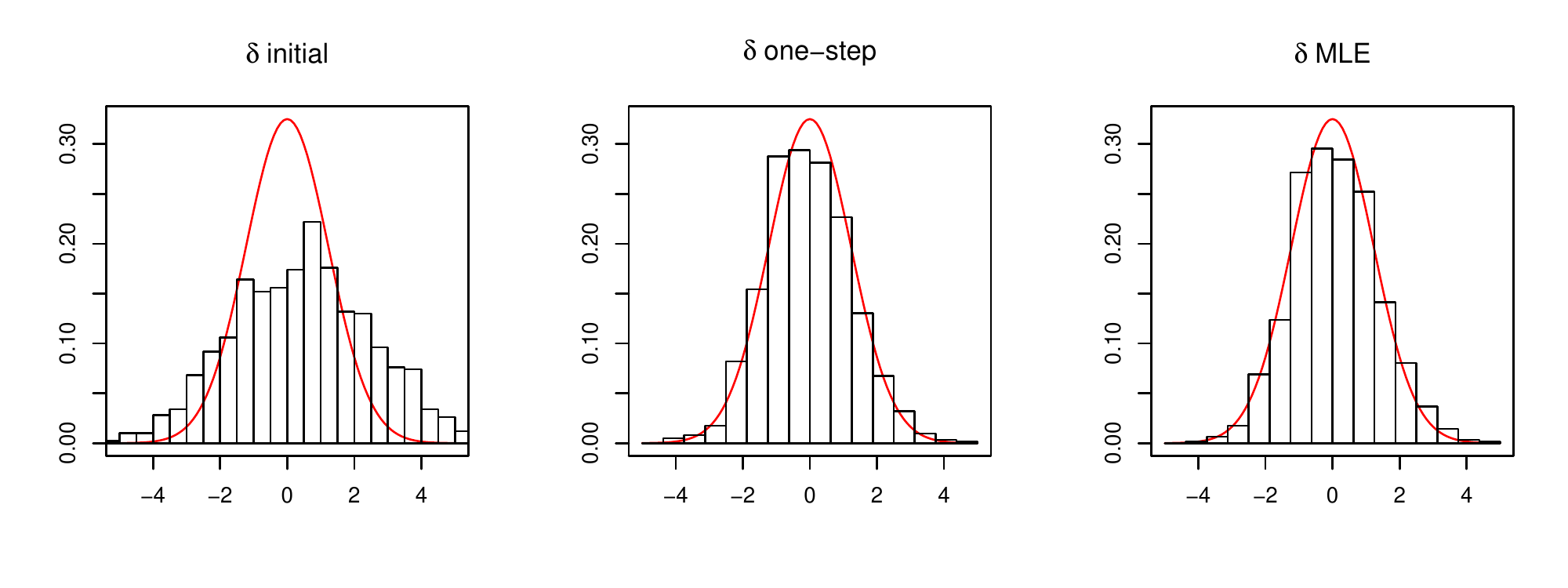}
  \end{center}
 \end{minipage}
 \begin{minipage}{0.49\hsize}
  \begin{center}
   \includegraphics[scale=0.38]{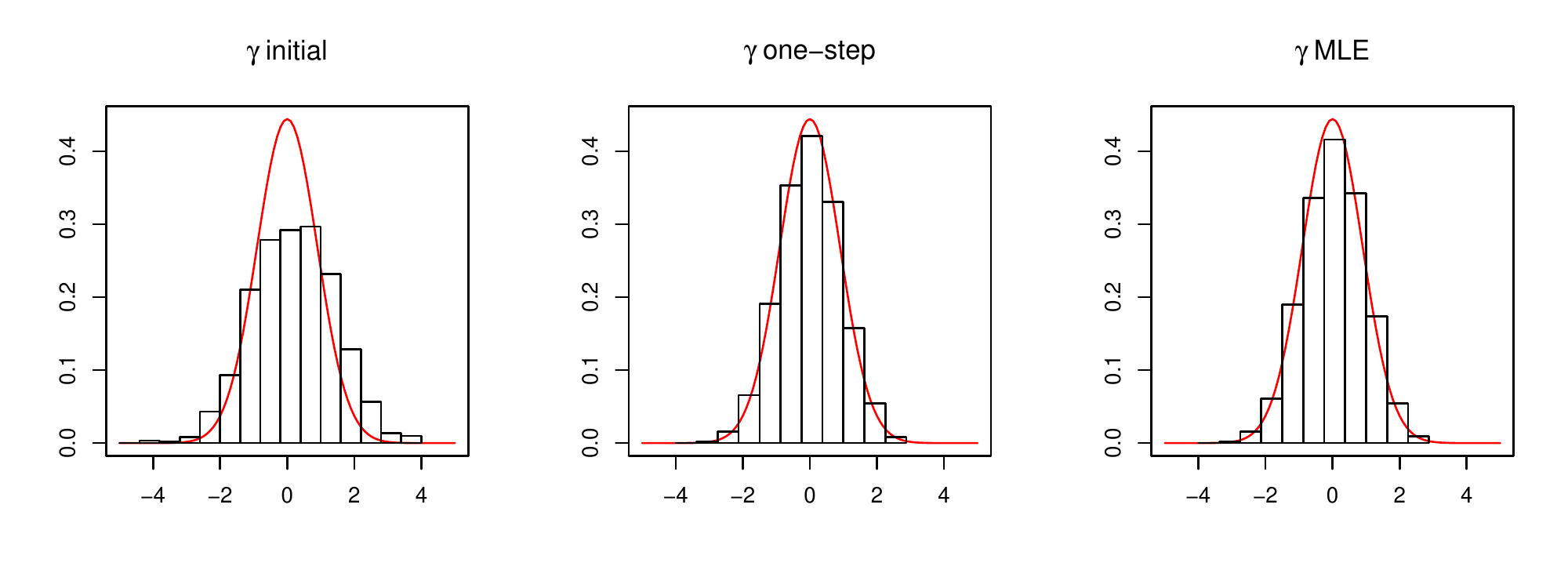}
  \end{center}
 \end{minipage}
 \\
   \caption{
   Case (i'): Histograms of the initial estimator $\tes^0$ (first and fourth columns), the one-step estimator $\tes^1$ (second and fifth columns), and the MLE $\tes$ (third and sixth columns). In each histogram panel, the solid red line shows the estimated asymptotically best possible normal distribution.
   }
\label{histograms_111}
\end{figure}

\begin{figure}[h]
 \centering
   \begin{minipage}{0.49\hsize}
  \begin{center}
   \includegraphics[scale=0.38]{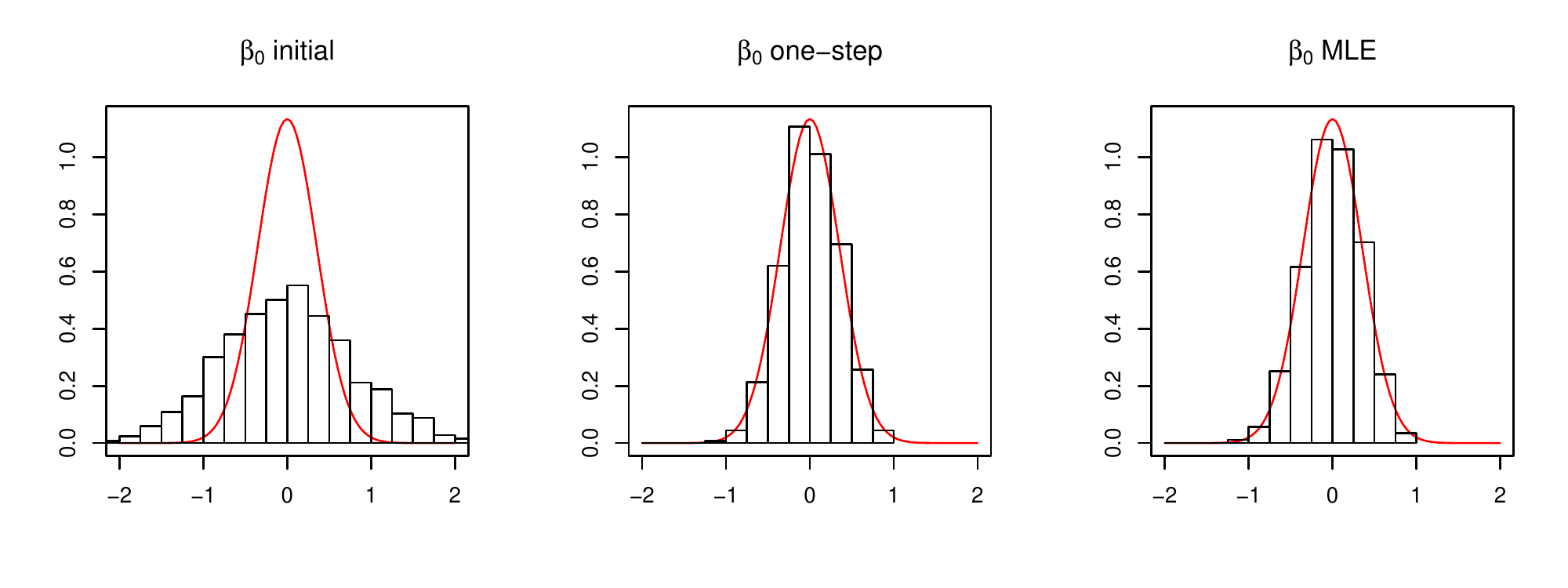}
  \end{center}
 \end{minipage}
 \begin{minipage}{0.49\hsize}
  \begin{center}
   \includegraphics[scale=0.38]{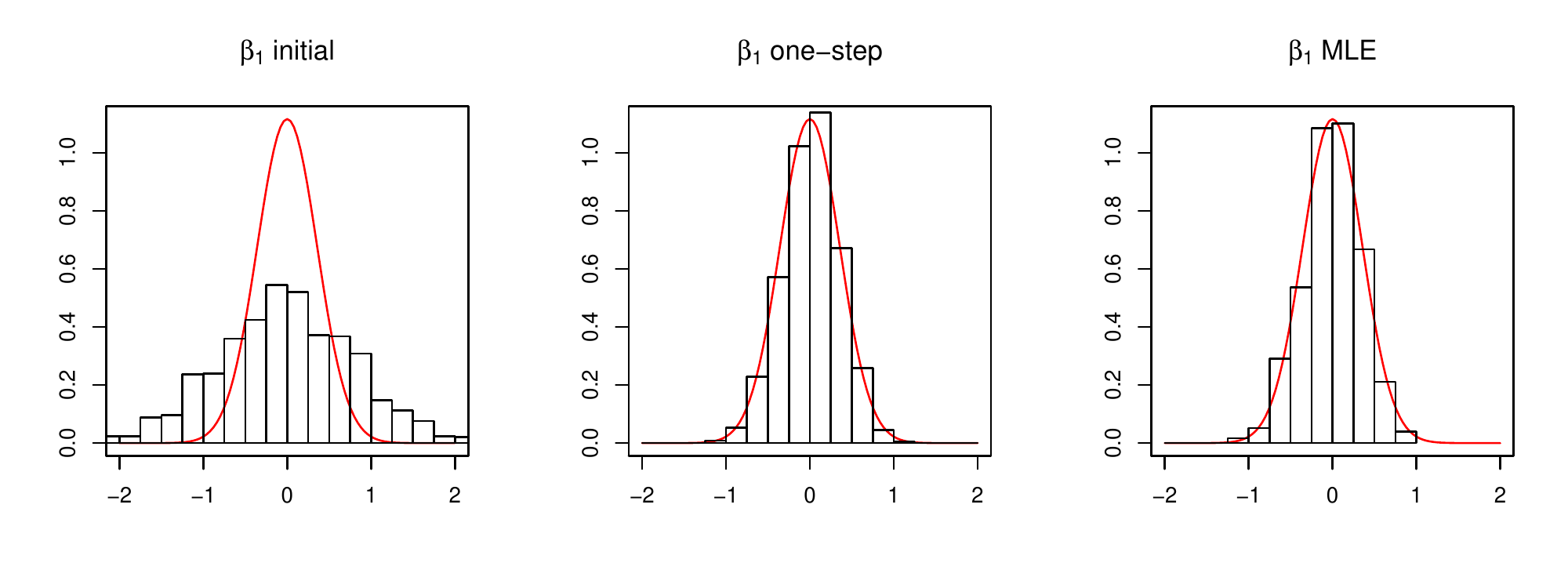}
  \end{center}
 \end{minipage}
 \\
  \begin{minipage}{0.49\hsize}
  \begin{center}
   \includegraphics[scale=0.38]{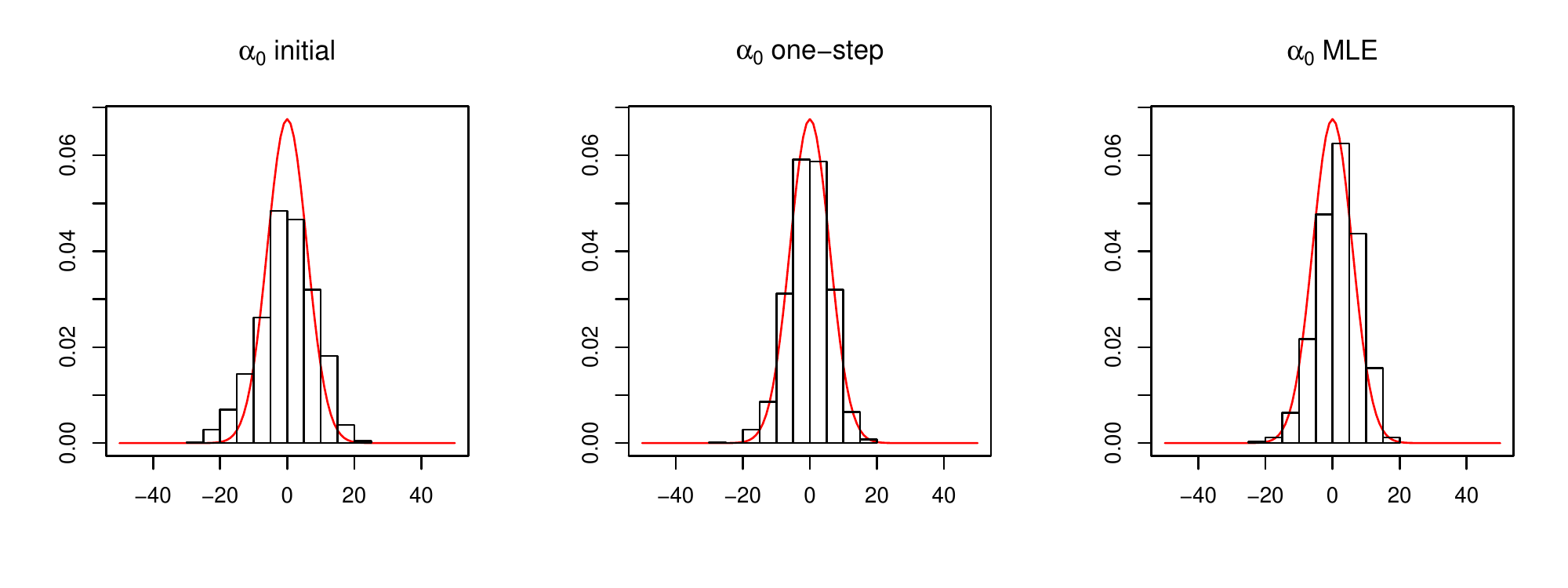}
  \end{center}
 \end{minipage}
 \begin{minipage}{0.49\hsize}
  \begin{center}
   \includegraphics[scale=0.38]{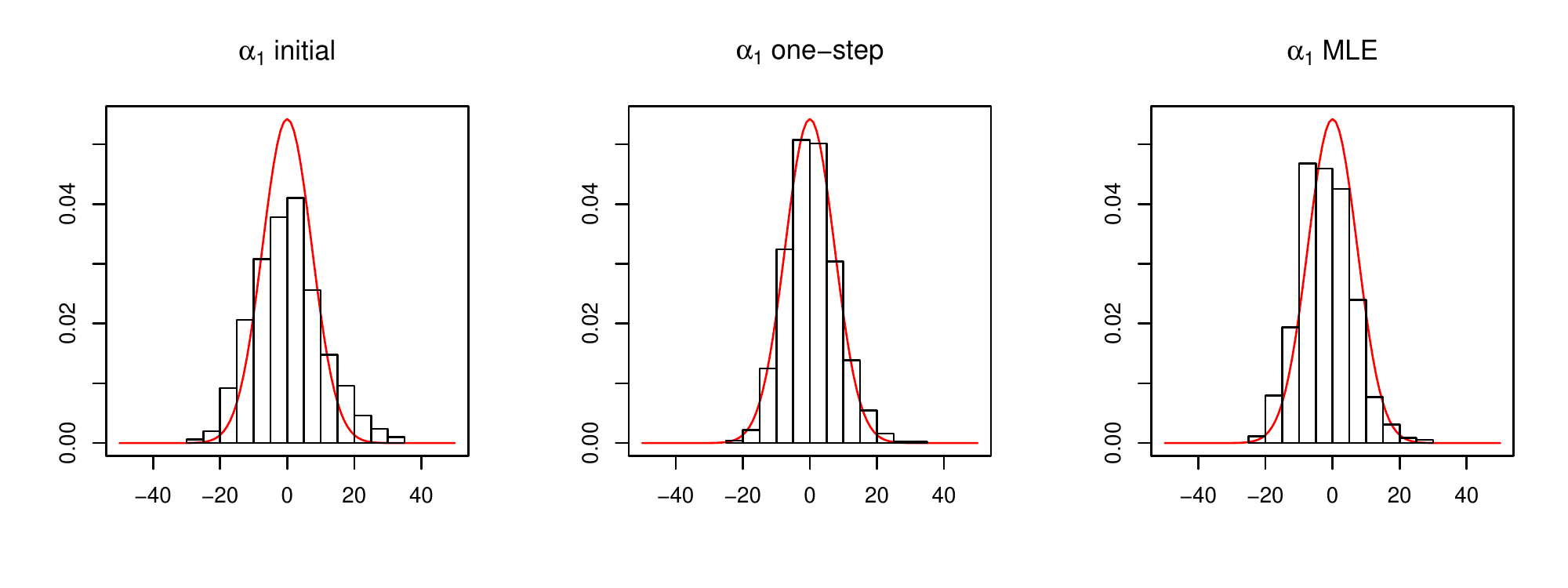}
  \end{center}
 \end{minipage}
 \\
 \begin{minipage}{0.49\hsize}
  \begin{center}
   \includegraphics[scale=0.38]{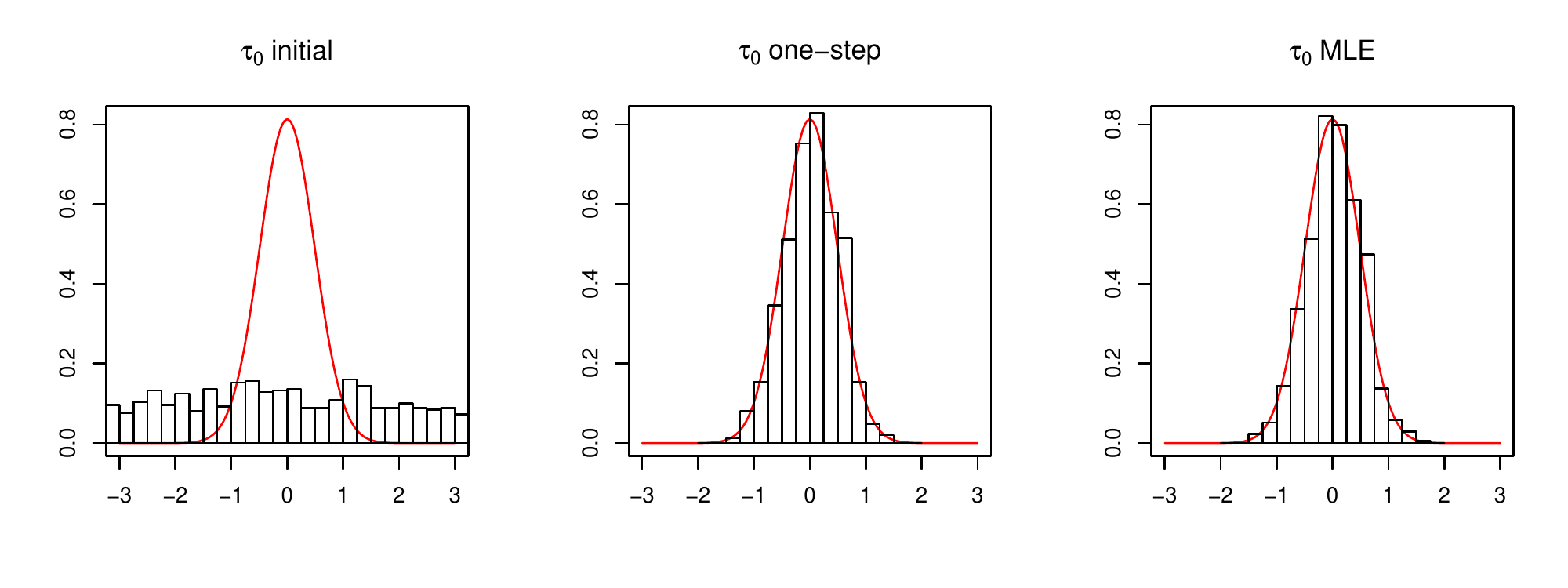}
  \end{center}
 \end{minipage}
 \begin{minipage}{0.49\hsize}
  \begin{center}
   \includegraphics[scale=0.38]{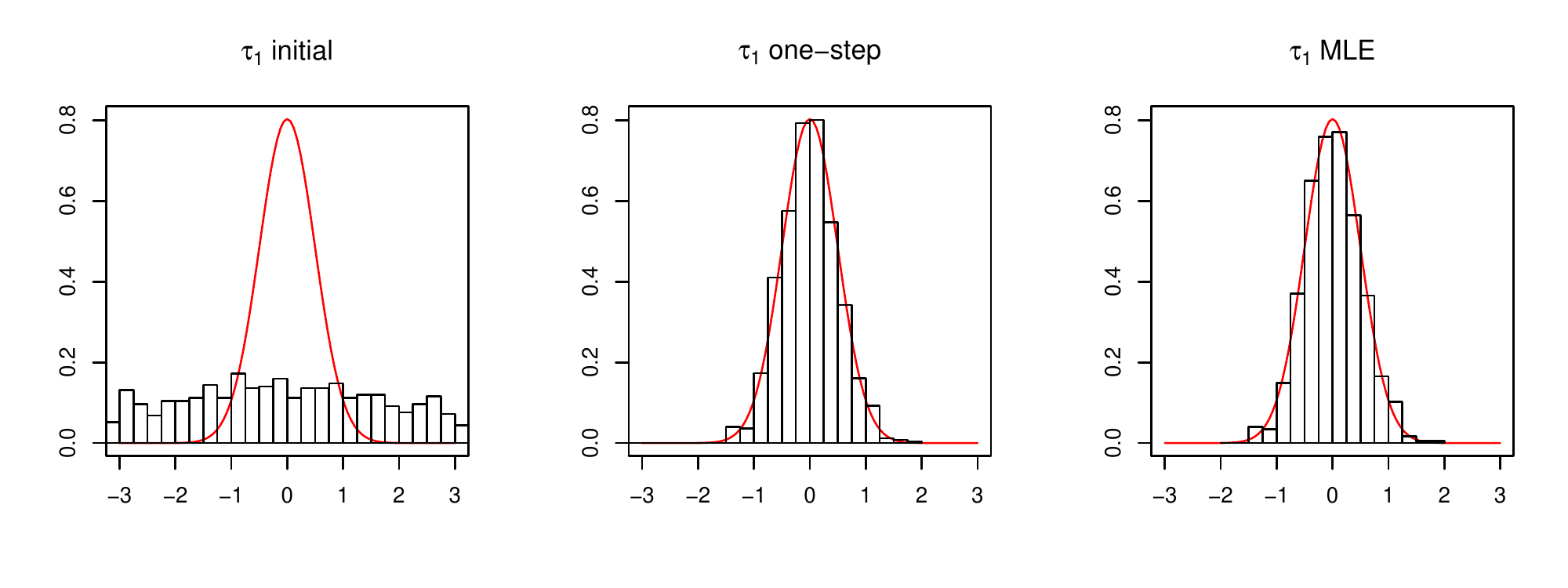}
  \end{center}
 \end{minipage}
 \\
  \begin{minipage}{0.49\hsize}
  \begin{center}
   \includegraphics[scale=0.38]{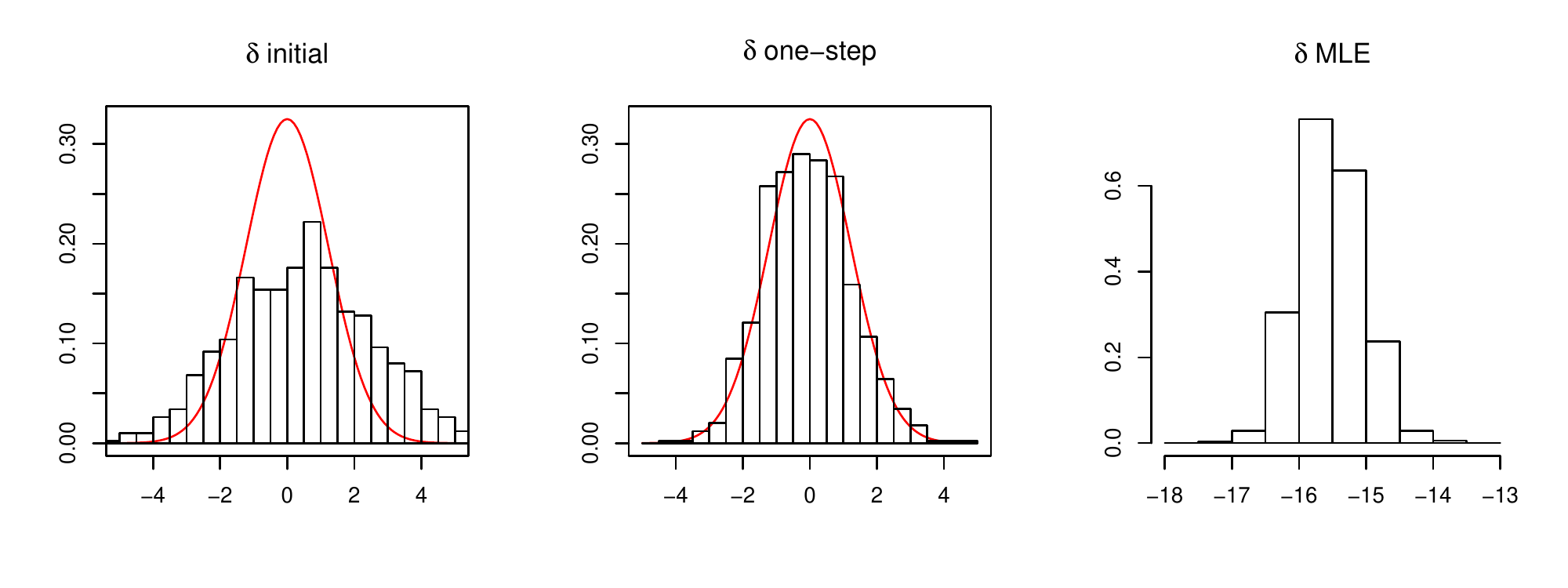}
  \end{center}
 \end{minipage}
 \begin{minipage}{0.49\hsize}
  \begin{center}
   \includegraphics[scale=0.38]{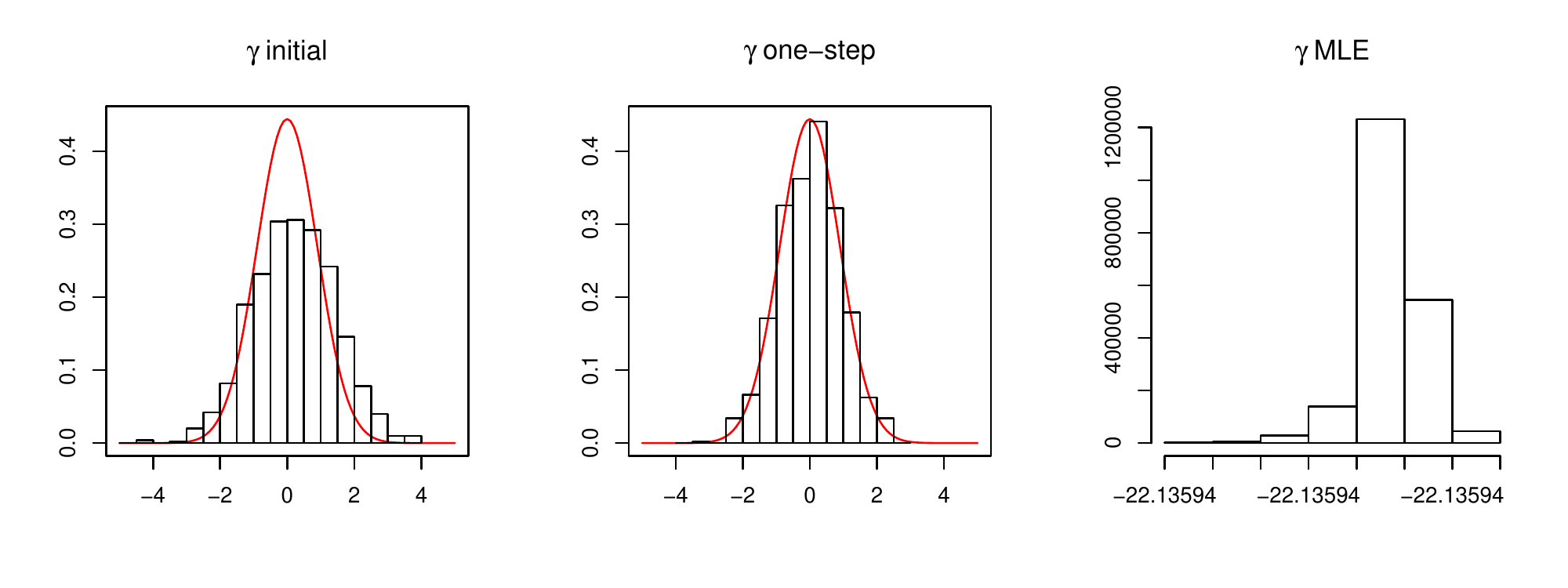}
  \end{center}
 \end{minipage}
 \\
   \caption{
   Case (ii'): 
     Histograms of the initial estimator $\tes^0$ (first and fourth columns), the one-step estimator $\tes^1$ (second and fifth columns), and the MLE $\tes$ (third and sixth columns). In each histogram panel, the solid red line shows the estimated asymptotically best possible normal distribution.
   }
     \label{histograms_222}
\end{figure}

In each case, we computed $\sqrt{N}(\hat{\xi}_N-\theta_0)$ for $\hat{\xi}_N = \tes^0$, $\tes^1$, and $\tes$, all being conducted $1000$-times Monte Carlo trials.
\rev{
To estimate $95\%$-coverage probabilities empirically as in Section \ref{sec_MLE.sim},
we computed the quantities $-\p_\theta^2\ell_N(\tes)$ and $-\p_\theta^2\ell_N(\tes^1)$ through the function $\theta \mapsto -\p_\theta^2\ell_N(\theta)$ for the approximately $95\%$-confidence intervals for each parameter.
The results are shown in Table \ref{coverage_MLE-ose}; therein, we obtained numerically unstable $4$ MLEs and $5$ one-step estimators for case (i') and $299$ MLEs and $6$ one-step estimators for case (ii'), and then computed the coverage probabilities based on the remaining cases.
In Figures \ref{histograms_111} and \ref{histograms_222} (for cases (i') and (ii'), respectively), we drew histograms of $\tes^1$ and $\tes$ together with those of the initial estimator $\tes^0$ for comparison.}
In each figure, the histograms in the first and fourth columns are those for $\tes^0$, those in the second and fifth columns for $\tes^1$, and those in the third and sixth columns for $\tes$, respectively;
the red solid line shows the zero-mean normal densities with the consistently estimated Fisher information for the variances.

\medskip

\begin{table}[h]
\centering
\begin{tabular}{cccccccccc}
  \hline
 & &$ \beta_0 $& $\beta_ 1$&$ \alpha_0 $& $\alpha_1 $&$ \tau_0 $& $\tau_1 $& $\del $& $\gam $\\
  \hline
Case (i') & $\tes$ & 0.960 & 0.959 & 0.944 & 0.942 & 0.957 & 0.947 & 0.948 & 0.943 \\ 
& $\tes^1$ & 0.960 & 0.958 & 0.943 & 0.943 & 0.952 & 0.949 & 0.943 & 0.947 \\[1mm]
Case (ii') & $\tes$ & 0.957 & 0.953 & 0.917 & 0.919 & 0.956 & 0.944 & 0& 0\\   
& $\tes^1$ & 0.960 & 0.959 & 0.943 & 0.943 & 0.952 & 0.949 & 0.943 & 0.947 \\ 
   \hline\\
\end{tabular}
\caption{The empirical $95\%$-coverage probabilities of the MLE and the one-step estimators in cases (i') and (ii') based on $1000$ trials;
MLE of $(\del,\gam)$ in case (ii') showed instability in numerical optimizations, \rev{while the one-step estimator is stable as in case (i').}}
\label{coverage_MLE-ose}
\end{table}

\medskip

\medskip

Here is a summary of the important findings.

\begin{itemize}

\item 
Approximate computation times for obtaining one set of estimates are as follows:
\begin{itemize}
\item[(i')] 0.2 seconds for $\tes^0$; \ 10 seconds for $\tes^1$; \ 2 minutes for $\tes$;
\item[(ii')] 0.2 seconds for $\tes^0$; \ 10 seconds for $\tes^1$; \ 9 minutes for $\tes$.
\end{itemize}
A considerable amount of reduction can be seen for $\tes^1$ compared with $\tes$.

\item 
About Figures \ref{histograms_111} and \ref{histograms_222}:

\begin{itemize}
\item In both cases (i') and (ii'), the inferior performance of $\tes^0$ is drastically improved by $\tes^1$, which in turn shows asymptotically equivalent behaviors to the MLE $\tes$.

\item On the one hand, as in Section \ref{sec_MLE.sim}, the MLE $\tes$ is much affected by the initial value for the numerical optimization, partly because of the non-convexity of the likelihood function $\ell_N(\theta)$; in Case (ii'), we observed the instability in computing the MLE of $(\del,\gam)$ (in the bottom panels in Figure \ref{histograms_222}), showing the local maxima problem.
On the other hand, we did not observe the local maxima problem in computing $\tes^0$ and the one-step estimator $\tes^1$ does not require an initial value for numerical optimization.
\end{itemize}

\end{itemize}

\medskip

In sum, $\tes^1$ is not only asymptotically equivalent to the efficient MLE but also much more robust in numerical optimization than the MLE.
It is recommended to use the one-step estimator $\tes^1$ against the MLE $\tes$ from both theoretical and computational points of view.

\medskip

\begin{figure}[h]
 \centering
  \begin{minipage}{0.49\hsize}
  \begin{center}
   \includegraphics[scale=0.38]{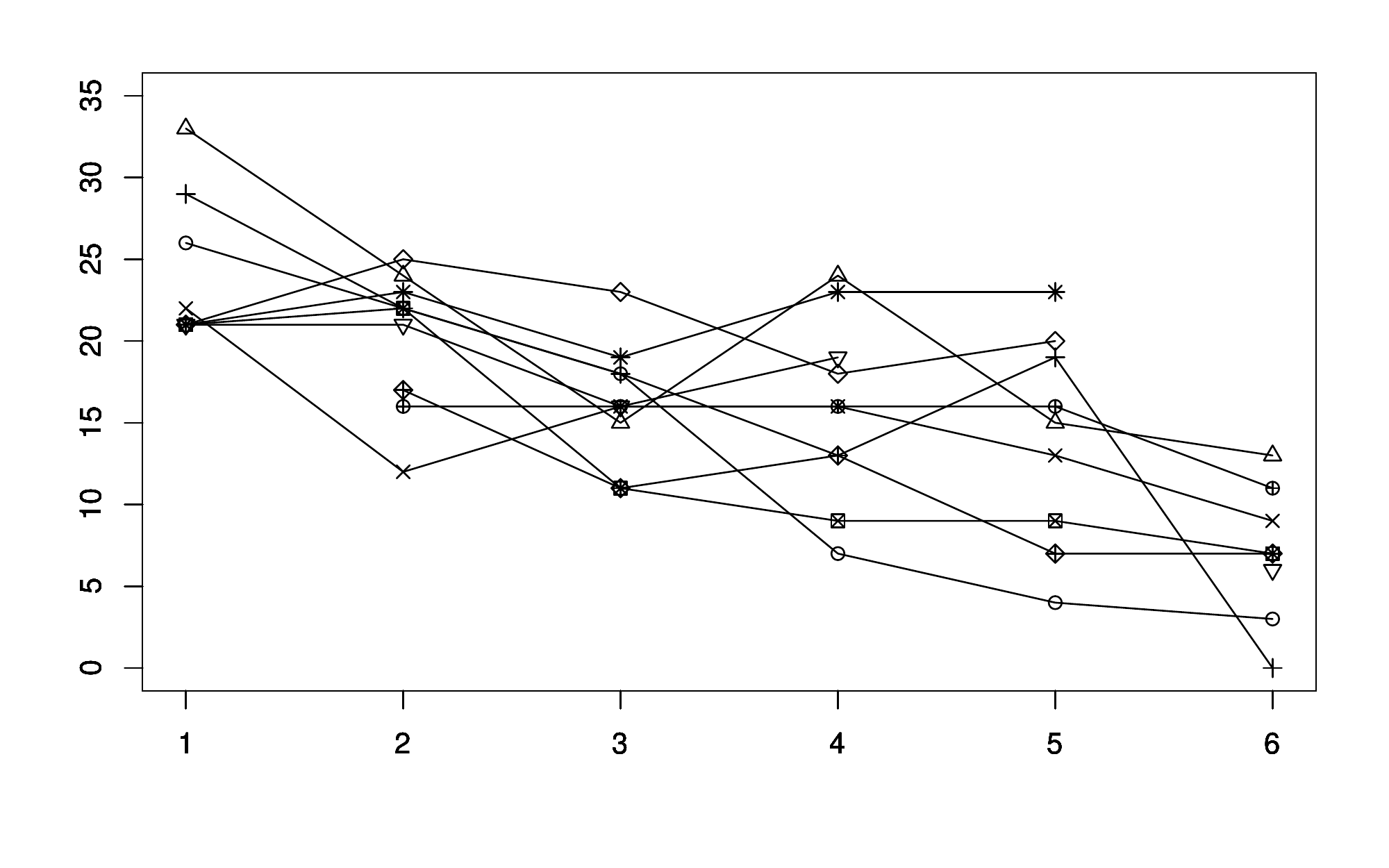}
  \end{center}
 \end{minipage}
 \begin{minipage}{0.49\hsize}
  \begin{center}
   \includegraphics[scale=0.38]{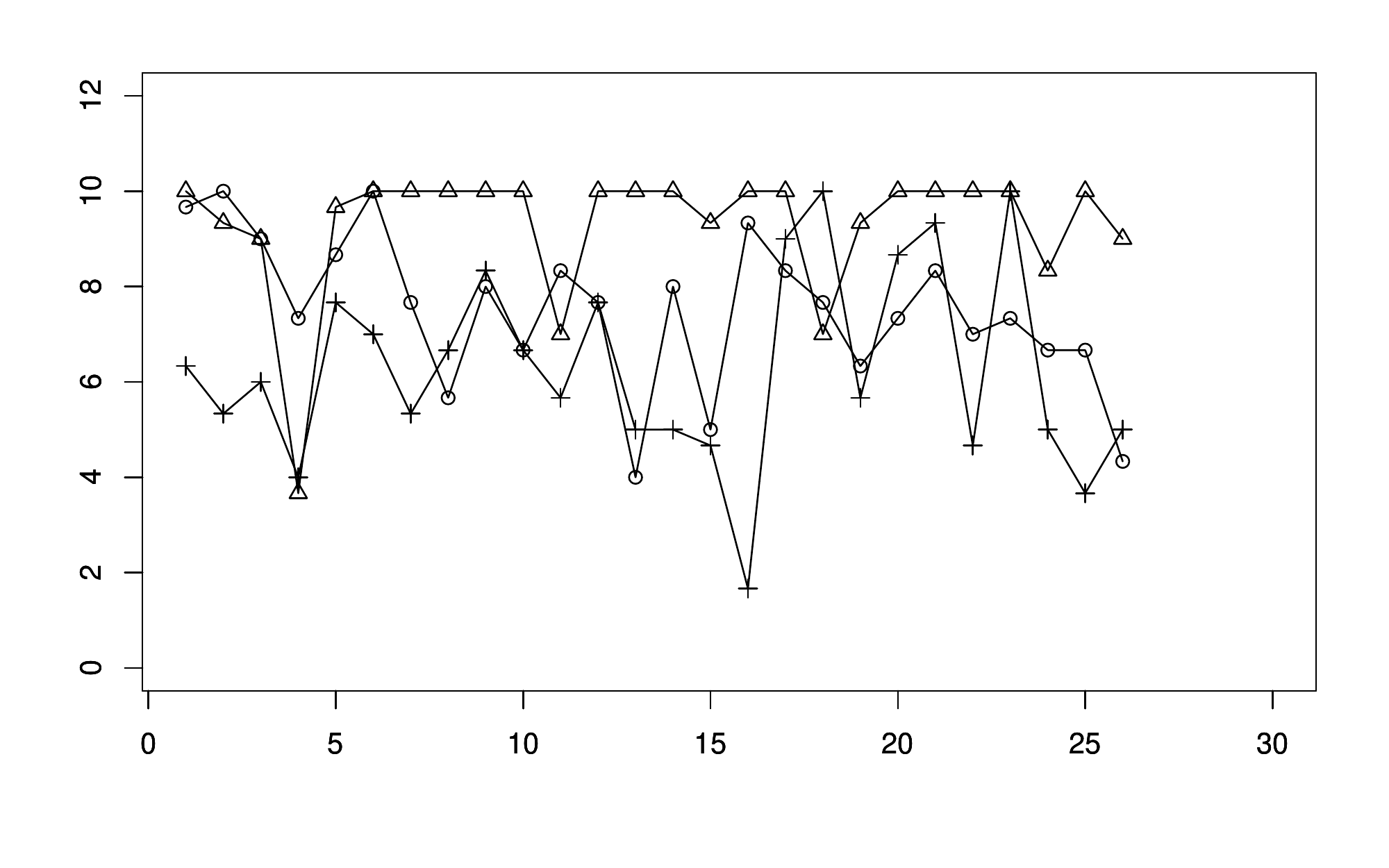}
  \end{center}
 \end{minipage}
 \\
   \caption{
   Data plots of \texttt{riesby\_example.dat} (left) and \texttt{posmood\_example.dat} (right) borrowed from the supplemental material of \cite{HedNor13}; the former shows data of $10$ patients over $6$ time points with a few missing values, and the latter does those of $3$ people over $26$ time points with no missing value.
   }
\label{dataplots_rda1}
\end{figure}

\begin{figure}[h]
 \centering
  \begin{minipage}{0.49\hsize}
  \begin{center}
   \includegraphics[scale=0.38]{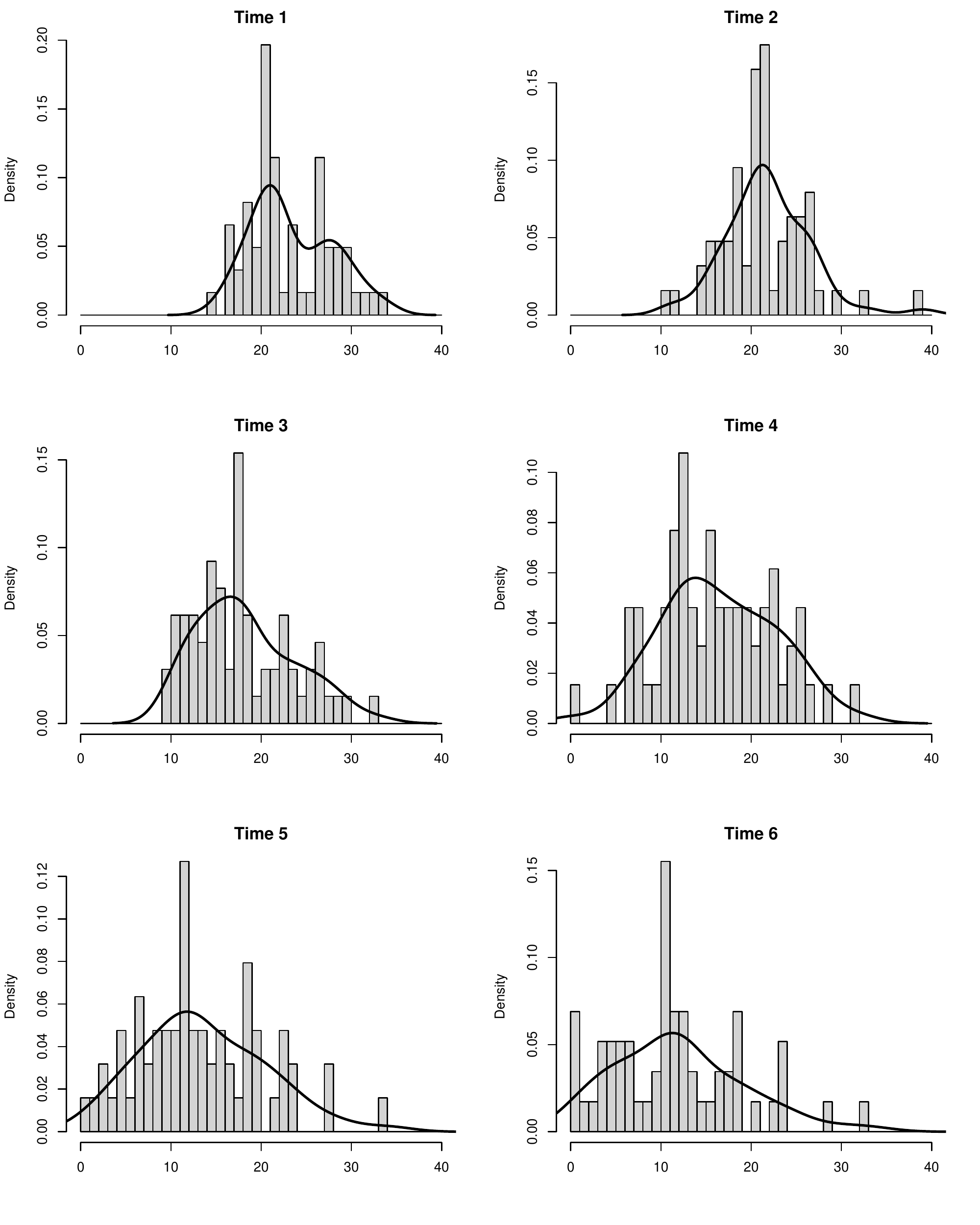}
  \end{center}
 \end{minipage}
 \begin{minipage}{0.49\hsize}
  \begin{center}
   \includegraphics[scale=0.38]{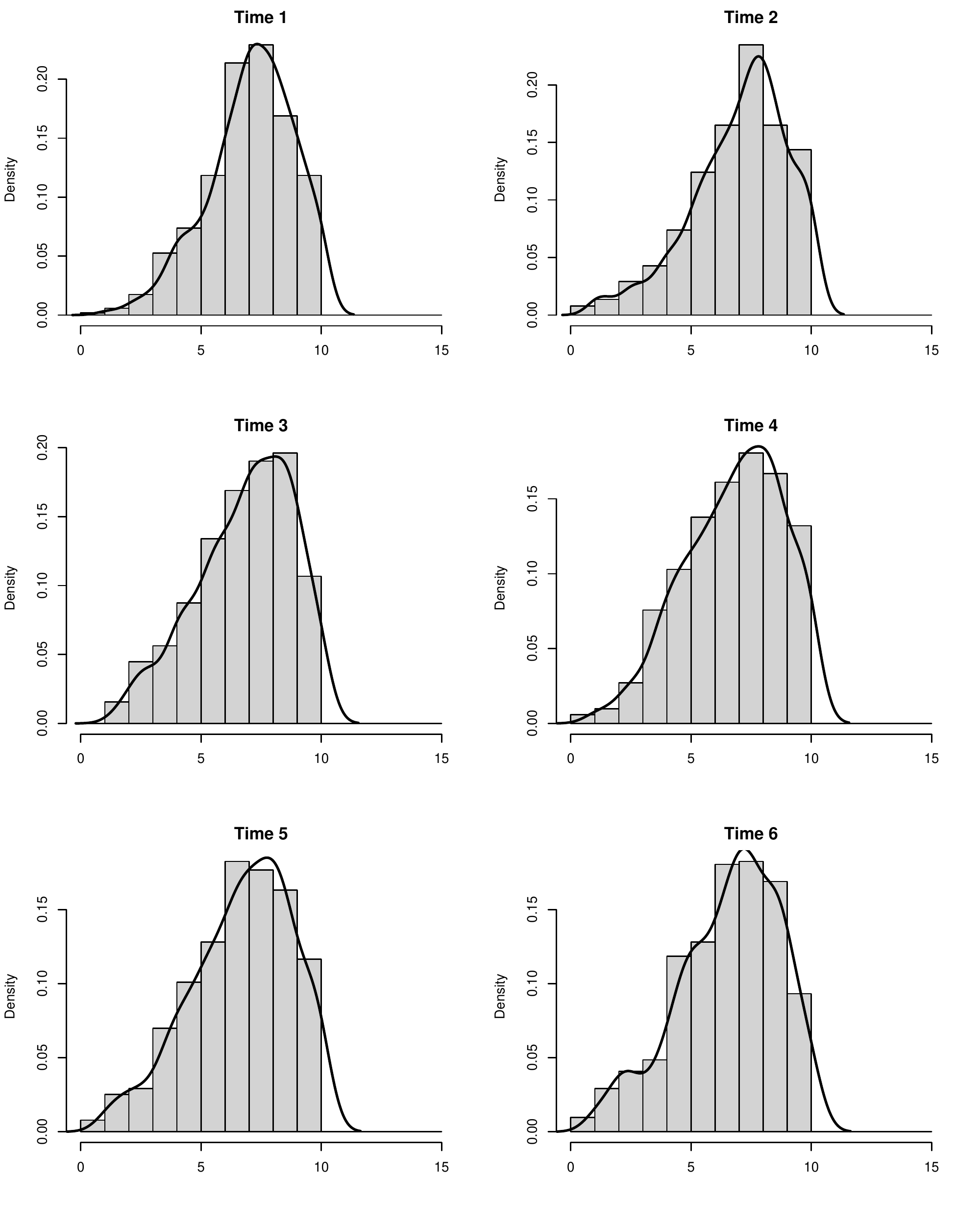}
  \end{center}
 \end{minipage}
 \\
   \caption{
   Histograms of data $(Y_{i1})_{i\le N},\dots,(Y_{i6})_{i\le N}$ for \texttt{riesby\_example.dat} (left) and \texttt{posmood\_example.dat} (right) borrowed from the supplemental material of \cite{HedNor13}.
   }
\label{dataplots_rda2}
\end{figure}

\rev{
We end this section with applications of the proposed one-step estimator $\tes^1$ for \eqref{model_ex.ose} to the two real data sets \texttt{riesby\_example.dat} and \texttt{posmood\_example.dat} borrowed from the supplemental material of \cite{HedNor13}.
Here are brief descriptions.
\begin{itemize}
\item \texttt{riesby\_example.dat} contains the Hamiltonian depression rating scale as $Y_{ij}$.
The covariates are given by 
$x_{ij}=(\texttt{intercept},\texttt{week},\texttt{edog})\in\mbbr\times\{0,1,2,\dots,5\}\times\{0,1\}$, 
$z_{ij}=(\texttt{intercept},\texttt{edog})$, and
$w_{ij}=(\texttt{intercept},\texttt{week})$.
Here, $N=66$ and the numbers of sampling times are $6$ with a few missing slots, and \texttt{edog} denotes the dummy variable for indicating whether the depression of the patient is endogenous ($=1$) or not ($=0$).
\item \texttt{posmood\_example.dat} contains the individual mood items as $Y_{ij}$; the items are pre-processed using factor analysis and take values $1$ to $10$ with higher ones indicating a higher level of positive mood.
The covariates are given by 
$x_{ij}=(\texttt{intercept},\texttt{alone},\texttt{genderf})\in\mbbr\times\{0,1\}\times\{0,1\}$, 
$z_{ij}=(\texttt{intercept},\texttt{alone})$, and
$w_{ij}=(\texttt{intercept},\texttt{alone})$.
Here, $N=515$ with no missing value, with approximately $34$ sampling times on average (ranging from $3$ to $58$). The variable \texttt{alone} and \texttt{genderf} respectively denote the dummy variables for indicating whether the person is alone ($=0$) or not ($=1$), which is time-varying, and whether the person is male ($=0$) or female ($=1$).
\end{itemize}
Figures \ref{dataplots_rda1} and \ref{dataplots_rda2} show some data plots and histograms, respectively; the former is positively skewed while the latter is negatively skewed.
We could apply our one-step estimation methods for these data sets, although they can be seen as categorical data (with a moderately large number of categories).
The results are given in Table \ref{table_rda}; the parameters $\beta_0$, $\al_0$, and $\tau_0$ denote the \texttt{intercept}.
The skewness mentioned above is reflected in the estimates of $\al_0$ and $\al_1$.
}

\medskip

\medskip

\medskip

\begin{table}[h]
\centering
\begin{tabular}{cccccccccc}
  \hline
\texttt{riesby\_example.dat} &$ \beta_0 $& $\beta_ 1$& $\beta_ 2$&$ \alpha_0 $& $\alpha_1 $&$ \tau_0 $& $\tau_1 $& $\del$ &$\gam$\\ 
  \cline{2-10}
 & 16.930 & -2.443 & -0.833 & 1.442 & 5.961 & -0.213 & 0.396 & 2.032 & 0.245 \\ 
\\[-2mm]
  \hline
\texttt{posmood\_example.dat} &$ \beta_0 $& $\beta_ 1$& $\beta_ 2$&$ \alpha_0 $& $\alpha_1 $&$ \tau_0 $& $\tau_1 $& $\del$ &$\gam$\\ 
\cline{2-10}
 & 6.608 & -0.229 & -0.295 & -0.125 & -0.016 & -0.332 & 0.093 & 4.517 & 0.606\\ 
   \hline\\[-2mm]
\end{tabular}
\caption{One-step estimates for the two data sets \texttt{riesby\_example.dat} and \texttt{posmood\_example.dat}.
It took 0.9 and 6.7 seconds, respectively.}
\label{table_rda}
\end{table}



\section{Concluding remarks}
\label{sec_cr}

%

We proposed a class of mixed-effects models with non-Gaussian marginal distributions which can incorporate random effects into the skewness and the scale simply and transparently through the normal variance-mean mixture.
The associated log-likelihood function is explicit and the MLE is asymptotically efficient (Remark \ref{rem_efficiency}) while computationally demanding and unstable. To bypass the numerical issue, we proposed the easy-to-use one-step estimator $\tes^1$, which turned out to not only attain a significant reduction of computation time compared with the MLE but also guarantee the asymptotic efficiency property.


Here are some remarks on important related issues.

\begin{enumerate}
\item \textit{Inter-individual dependence structure.}
A drawback of the model \eqref{gh-melsr} is that its inter-individual dependence structure is not flexible enough.
Specifically, let us again note 
the following covariance structure for $j,k\le n_i$:
\begin{equation}
\cov_{\theta}[Y_{ij}, Y_{ik}] = s_{ij}(\al) s_{ik}(\al) \var_{\theta}[v_i] = c(\theta') s_{ij}(\al) s_{ik}(\al).
\nonumber
\end{equation}
This in particular implies that $Y_{i1},\dots,Y_{i n_i}$ cannot be correlated as long as $s(z,\al)\equiv 0$.
Nevertheless, it is formally straightforward to extend the model \eqref{gh-melsr} so that the distributional structure of $Y_i \in\mbbr^{n_i}$ obeys the multivariate GH distribution for each $\mcl(Y_i)$ with a \textit{non-diagonal} scale matrix.
To mention it briefly, suppose that the vector of a sample $Y_i=(Y_{i1},\dots,Y_{i n_i})\in \mbbr^{n_i}$ from $i$th individual is given by the form
\begin{eqnarray}
Y_i = x_i\beta+ s(z_i,\alpha)v_i + \Lambda(w_i,\tau)^{1/2} \sqrt{v_i}\,\epsilon_{i}.
\nonumber
\end{eqnarray}
Here $v_1,\ldots,~v_N\sim \text{i.i.d.}~GIG(\lambda,\delta,\gamma)$ as before, while we now incorporated the scale matrix $\Lam(w_i,\tau)$ which should be positive definite and symmetric, but may be \textit{non-diagonal}. Then, the dependence structure of $Y_{i1},\dots,Y_{i n_i}$ can be much more flexible than \eqref{gh-melsr}.

%

\item \textit{Forecasting random-effect parameters.}
In the familiar Gaussian linear mixed-effects model of the form $Y_i=X_i\beta+Z_i b_i + \ep_i$, the empirical Bayes predictor of $v_i$ is given by $\hat{b}_i := \E_\theta[b_i|Y_i]|_{\theta=\tes}$.
One of the analytical merits of our NVMM framework is that the conditional distribution $\mcl(v_i|Y_i=y_i)$ of $v_i$ is given by $GIG(\nu_i,\eta_i,\psi_i)$, where
\begin{align}
\nu_i = \nu_i(\theta) &:= \lambda-\frac{n_i}{2}, \nn\\
\eta_i = \eta_i(\theta) &:= \sqrt{\delta^2+(y_i - x_i\beta)^{\top}\Lambda(w_,\tau)^{-1}(y_i-x_i\beta)}, \nn\\
\psi_i = \psi_i(\theta) &:= \sqrt{\gamma^2 + s_i(\al)^\top \Lambda(w_i,\tau)^{-1}s_i(\al)}.
\nonumber
\end{align}
This is a direct consequence of the general results about the multivariate GH distribution; see \cite{EbeHam04} and the references therein for details.
As in the Gaussian case mentioned above, we can make use of
\begin{equation}
\hat{v}_i := \E_\theta[v_i|Y_i=y_i]|_{\theta=\tes} 
=\frac{K_{\hat{\nu}_i+1}(\hat{\eta}_i\hat{\psi}_i)}{K_{\hat{\nu}_i}(\hat{\eta}_i\hat{\psi}_i)} \frac{\hat{\eta}_i}{\hat{\psi}_i},
\nonumber
\end{equation}
where $\hat{\nu}_i := \nu_i(\tes)$, $\hat{\eta}_i := \eta_i(\tes)$, and $\hat{\psi}_i := \psi_i(\tes)$;
formally $\tes$ could be replaced by the one-step estimator $\tes^1$. Then, it would be natural to regard
\begin{equation}
\hat{Y}_{ij} := x_{ij}'\bes + s(z_{ij}',\aes) \hat{v}_i
\nonumber
\end{equation}
as a prediction value of $Y_{ij}$ at $(x_{ij}',z_{ij}')$.
This includes forecasting the value of $i$th individual at a future time point.

\rev{
\item \textit{Lack of fit and model selection.}
In relation to Remark \ref{rem_misspecification}, based on the obtained asymptotic-normality results, we can proceed with lack-of-fit tests, such as the likelihood-ratio test, the score test, and the Wald test;
typical forms are $s(z,\al)=\sum_{l=1}^{p_\al}\al_l s_l(z)$ and $\sig(w,\tau)=\exp\{\sum_{m=1}^{p_\tau}\tau_m \sig_m(w)\}$, with given basis functions $s_l(z)$ and $\sig_m(w)$. In that case, we can estimate $p$-value for each component of $\theta$, say, by $2\Phi(-|\hat{B}_{k,N} \hat{\theta}_{k,N}|)$ for $\theta_k$ where $\hat{B}_{k,N}:=[(-\p_{\theta}^2\ell_N(\tes))^{-1}]_{kk}^{-1/2}$.
Alternatively, one may consider information criteria such as the conditional AIC \cite{VaiBla05} and the BIC-type one \cite{DelLavPou14}.
To develop these devices in rigorous ways, we will need to derive several further analytical results: the uniform integrability of $(\|\sqrt{N}(\tes -\tz)\|^2)_n$ for the AIC, the stochastic expansion for the marginal likelihood function for the BIC, and so on.
}

\end{enumerate}

\bigskip

\noindent
{\bf Acknowledgement.}
The authors should like to thank the editors and the anonymous reviewers for their valuable comments, which led to substantial improvement of the paper.
This work was partly supported by JST CREST Grant Number JPMJCR2115, and by JSPS KAKENHI Grant Number 22H01139, Japan (HM).

\medskip

\noindent
{\bf Conflict of interest.}
The author declares that there is no conflict of interest.

\bigskip


\def\cprime{$'$} \def\polhk#1{\setbox0=\hbox{#1}{\ooalign{\hidewidth
  \lower1.5ex\hbox{`}\hidewidth\crcr\unhbox0}}} \def\cprime{$'$}
  \def\cprime{$'$}

\appendix

\section{GIG and GH distributions}
\label{sec_gig&gh}

Let $K_\lam(t)$ denote the modified Bessel function of the second kind ($\nu\in\mbbr$, $t>0$):
\begin{equation}
K_{\nu}(t)=\frac{1}{2}\int_0^{\infty}s^{\nu-1}\exp\left\{-\frac{t}{2}\left(s+\frac{1}{s}\right)\right\}ds.
\nonumber
\end{equation}
We have the following recurrence formulae \cite{AbrSte92}:
$K_{\nu+1}(t)=\frac{2\nu}{t}K_{\nu}(t)+K_{\nu-1}(t)$ and $K_{\nu-1}(t)+K_{\nu+1}(t)=-2\partial_t K_{\nu}(t)$.
It follows that $K_{\nu}(t)$ is monotonically decreasing and that $\partial_t K_{\nu}(t)=-K_{\nu-1}(t)-\frac{\nu}{t}K_{\nu}(t)$.
Further, we have
\begin{align}
\partial_t \log K_{\nu}(t) &= - \frac{K_{\nu-1}(t)}{K_{\nu}(t)} -\frac{\nu}{t} =: -R_\nu(t) -\frac{\nu}{t},
\label{def_R}\\
\partial^2_t \log K_{\nu}(t) &= - \frac{1}{K^2_{\nu}(t)}\left(K^2_{\nu-1}(t)-K_{\nu-2}(t)K_{\nu}(t)\right) -\frac{1}{t}\frac{K_{\nu-1}(t)}{K_{\nu}(t)}+\frac{\nu}{t^2}
\nn\\
&=: -S_\nu(t) - \frac{1}{t}R_\nu(t) + \frac{\nu}{t^2}.
\label{def_S}
\end{align}
The following asymptotic behavior holds:
\begin{equation}
K_{\nu}(t)=\sqrt{\frac{\pi}{2t}}\exp(-t)\{1+(4\nu^2-1)O(t^{-1})\},\qquad t\to\infty.
\nonumber
\end{equation}

The generalized inverse Gaussian (GIG) distribution $GIG(\lam,\del,\gam)$ on $\mbbrp$ is defined by the density:
\begin{equation}
p_{GIG}(z;\lambda,\delta,\gamma)=\frac{(\gamma/\delta)^{\lambda}}
{2K_{\lambda}(\gamma\delta)} z^{\lambda -1}
\exp\left\{-\frac{1}{2}\left(\frac{\delta^{2}}{z}+\gamma^{2}z\right)\right\},\qquad z>0.
\label{gig-d}
\end{equation}
The region of admissible parameters is given by the union of 
$\{(\lam,\del,\gam):\,\lam>0,\,\delta\ge 0,\, \gamma >0\}$, 
$\{(\lam,\del,\gam):\,\lam=0,\,\delta> 0,\, \gamma >0\}$, and
$\{(\lam,\del,\gam):\,\lam<0,\,\delta>0,\, \gamma >0\}$, 
according to the integrability of $p_{GIG}$ at the origin and $+\infty$.

The generalized hyperbolic (GH) distribution denoted by $GH(\lam,\al,\beta,\del,\mu)$ is defined as the distribution of the normal variance-mean mixture $Y$ with respect to $Z \sim GIG(\lam,\del,\gam)$:
\begin{equation}
Y=\mu+\beta Z+\sqrt{Z}\eta,
\nn
\end{equation}
where $\al:=\sqrt{\beta^2+\gam^2}$ 
and $\eta\sim N(0,1)$ independent of $Z$.
By the conditional Gaussianity $\mcl(Y|Z=z)=N(\mu+\beta z, z)$, the density is calculated as follows:
\begin{align}
& p_{GH}(y;\lambda,\alpha,\beta,\delta,\mu) \nn\\
&=
\int_{0}^{\infty} \frac{1}{\sqrt{2\pi z}}
\exp\left(-\frac{1}{2z}(y-\mu-\beta z)^2\right) p_{GIG}(z;\lambda,\delta,\gamma)dz
\nn\\
&=\frac{\left(\alpha^2-\beta^2\right)^{\lambda/2}\sqrt{\delta^2+(y-\mu)^2}^{\lambda-1/2}}{\sqrt{2\pi}\alpha^{\lambda-1/2}\delta^\lambda K_\lambda\left(\delta\sqrt{\alpha^2-\beta^2}\right)}K_{\lambda-\frac{1}{2}}\left(\alpha\sqrt{\delta^2+(y-\mu)^2}\right)\exp[\beta(y-\mu)].
\nn
\end{align}
The region of admissible parameters is given by the union of 
$\{(\lam,\al,\beta,\del,\mu):\,\lam>0,\,\delta\ge 0,\, \al>|\beta|\}$, 
$\{(\lam,\al,\beta,\del,\mu):\,\lam=0,\,\delta> 0,\, \al>|\beta|\}$, and
$\{(\lam,\al,\beta,\del,\mu):\,\lam<0,\,\delta>0,\, \al\ge|\beta|\}$.
The mean and variance of $Y\sim GH(\lam,\al,\beta,\del,\mu)$ are given by
\begin{align}
\E[Y] &= \mu+\frac{\delta\beta K_{\lambda+1}(\delta\gamma)}{\gamma K_\lambda  (\delta\gamma)}, \nn\\
\var[Y] &= \frac{\delta K_{\lambda+1}(\delta\gamma)}{\gamma K_\lambda  (\delta\gamma)}+\frac{\beta^2\delta^2}{\gamma^2}\left[   \frac{ K_{\lambda+2}(\delta\gamma)}{ K_\lambda  (\delta\gamma)}   -\left(\frac{ K_{\lambda+1}(\delta\gamma)}{K_\lambda  (\delta\gamma)}\right) ^2 \right].
\nonumber
\end{align}
See \cite{EbeHam04} for further details of the GIG and GH distributions.

The normal inverse Gaussian (NIG) distribution is one of the popular subclasses of the GH-distribution family:
$NIG(\al,\beta,\del,\mu):=GH(-1/2,\al,\beta,\del,\mu)$, where $GIG(-1/2,\del,\gam)$ corresponds to the inverse Gaussian distribution.
The $NIG(\al,\beta,\del,\mu)$-density is given by
\begin{equation}
p_{NIG}(x;\alpha,\beta,\delta,\mu)=\frac{\alpha\delta}{\pi}\exp\left( 
\delta\sqrt{\alpha^2-\beta^2}+\beta(x-\mu)\right) \frac{K_1\left(\alpha\sqrt{\delta^2+(x-\mu)^2}\right)}{\sqrt{\delta^2+(x-\mu)^2}}.
\nonumber
\end{equation}
All of the mean $M$, variance $V$, skewness $S$, and kurtosis $K$ of $NIG(\al,\beta,\del,\mu)$ are explicitly given:
\begin{eqnarray}
M=\mu+\frac{\beta\delta}{(\alpha^2-\beta^2)^{\frac{1}{2}}},
\quad V=\frac{\delta\alpha^2}{(\alpha^2-\beta^2)^{\frac{3}{2}}},
\quad S=\frac{3\beta}{\alpha\sqrt{\delta}(\alpha^2-\beta^2)^{\frac{1}{4}}},
\quad K=\frac{3\alpha^2+4\beta^2}{\alpha^2\delta(\alpha^2-\beta^2)^{\frac{1}{2}}}.
\nonumber
\end{eqnarray}
Inverting these expressions gives
\begin{eqnarray}
\gamma=\frac{3}{\sqrt{V}\sqrt{3K-5S^2}},
\quad \beta=\frac{S\sqrt{V}\gamma^2}{3},
\quad \alpha=\sqrt{\gamma^2+\beta^2},
\quad \delta=\frac{V\gamma^3}{\gamma^2+\beta^2},
\quad \mu=M-\frac{\beta\delta}{\gamma},
\nn
\end{eqnarray}
from which one can consider the method-of-moments estimation of $(\al,\beta,\del,\mu)$ based on the empirical counterparts of $M$, $V$, $S$, and $K$. One should note that the empirical quantity $3\hat{K}_n-5\hat{S}_n^2$ has to be positive, which may fail in a finite sample and for such a data set the MLE would be also non-computable or unstable.
In \cite{YooKimSon20}, the estimation problem for the i.i.d. NIG model was studied from the computational point of view;
the paper also introduced the change of variables for the parameters to sidestep the positivity restriction, resulting in stabilized results in numerical experiments.

\section{Likelihood function}
\label{sec_log-LF}

\subsection{Derivation}
\label{sec_log-LF-derivation}

Writing $\theta_1=(\beta,\alpha,\tau)$ and $\theta_2=(\lambda,\delta,\gamma)$, and using the obvious notation, we obtain
\begin{align}
\ell_n(\theta)&=\log p_{\theta}(Y_1,\ldots,Y_n)\nonumber\\
&=\log\int\cdots\int p_{\theta_1}(Y_1,Y_2,\ldots,Y_n|v_1,\ldots,v_n)\prod_{i=1}^N p_{\theta_2} (v_i)dv_i\nonumber\\
&=\log\left(\int\dots\int \prod_{i=1}^N p_{\theta_1} (Y_i|v_i)\prod_{i=1}^N p_{\theta_2} (v_i)dv_i\right)
\nonumber\\
&=\sum_{i=1}^N   \log\left[\int \Bigg(\prod_{j=1}^{n_i} p_{\theta_1} (Y_{ij}|v_i)\Bigg)p_{\theta_2} (v_i)dv_i\right]\nonumber\\
&= \sum_{i=1}^N   \log\left[\int \left( \prod_{j=1}^{n_i} \frac{1}{\sqrt{2\pi\sigma^2_{ij}(\tau)}}v_i^{-\frac{1}{2}}\exp\left[-\frac{1}{2\sigma^2_{ij}(\tau)}(Y_{ij}-x^{\top}_{ij}\beta-  s_{ij}(\alpha)v_i)^2\right]\right)p_{\theta_2}(v_i) dv_i\right]\nonumber\\
&= \sum_{i=1}^N    \log\Biggl[ (2\pi)^{-\frac{n_i}{2}}\Bigg(\prod_{j=1}^{n_i} \sigma^2_{ij}(\tau)\Bigg)^{-\frac{1}{2}}\int_0^{\infty} v_i^{-\frac{n_i}{2}}\nonumber\\
&{}\qquad\times\prod_{j=1}^{n_i}\exp\left[-\frac{1}{2v_i\sigma^2_{ij}(\tau)}\left\{(Y_{ij}-x^{\top}_{ij}\beta)^2+s^2_{ij}(\alpha)v^2_i-2  s_{ij}(\alpha)v_i(Y_{ij}-x^{\top}_{ij}\beta)\right\}\right]p_{\theta_2}(v_i) dv_i\Biggr]\nonumber\\
&= \sum_{i=1}^N \log\Biggl[ (2\pi)^{-\frac{n_i}{2}}\Bigg(\prod_{j=1}^{n_i} \sigma^2_{ij}(\tau)\Bigg)^{-\frac{1}{2}}\int_0^{\infty} v_i^{-\frac{n_i}{2}}\nonumber\\
&{}\qquad \times\prod_{j=1}^{n_i}\exp\left[-\frac{1}{2v_i\sigma^2_{ij}(\tau)}\left\{(Y_{ij}-x^{\top}_{ij}\beta)^2+s^2_{ij}(\alpha)v^2_i\right\}+ \frac{ s_{ij}(\alpha)}{\sigma^2_{ij}(\tau)}(Y_{ij}-x^{\top}_{ij}\beta)\right]p_{\theta_2}(v_i) dv_i\Biggr]\nonumber\\
&= \sum_{i=1}^N    \log\Biggl[ (2\pi)^{-\frac{n_i}{2}}\Bigg(\prod_{j=1}^{n_i} \sigma^2_{ij}(\tau)\Bigg)^{-\frac{1}{2}}\prod_{j=1}^{n_i}\exp\left(\frac{  s_{ij}(\alpha)}{\sigma^2_{ij}(\tau)}(Y_{ij}-x^{\top}_{ij}\beta)\right)\nonumber\\
&{}\qquad\times \int_0^{\infty} v_i^{-\frac{n_i}{2}}\prod_{j=1}^{n_i}\exp\left\{-\frac{1}{2\sigma^2_{ij}(\tau)}\left(\frac{(Y_{ij}-x^{\top}_{ij}\beta)^2}{v_i}+s^2_{ij}(\alpha)v_i\right)\right\} p_{\theta_2}(v_i) dv_i\Biggr]\nonumber\\
&= \sum_{i=1}^N \log\Biggl[ 
(2\pi)^{-\frac{n_i}{2}}\Bigg(\prod_{j=1}^{n_i} \sigma^2_{ij}(\tau)\Bigg)^{-\frac{1}{2}}\exp\left( \sum_{j=1}^{n_i}\frac{ s_{ij}(\alpha)(Y_{ij}-x^{\top}_{ij}\beta)}{\sigma^2_{ij}(\tau)}\right) 
\nonumber\\
&{}\qquad \times\int_0^{\infty} v_i^{-\frac{n_i}{2}}
\exp\left\{ -\frac{1}{2}\left(\sum_{j=1}^{n_i}\frac{(Y_{ij}-x^{\top}_{ij}\beta)^2}{v_i\sigma^2_{ij}(\tau)}+\sum_{j=1}^{n_i}\frac{s^2_{ij}(\alpha)}{\sigma^2_{ij}(\tau)}v_i\right) \right\} p_{\theta_2}(v_i) dv_i\Biggr]\nonumber\\
&= \sum_{i=1}^N \log\Bigg[
C_i(\al,\beta,\tau) 
\int_0^{\infty} v_i^{-\frac{n_i}{2}}\exp\Bigg\{-\frac{1}{2}\Bigg(\sum_{j=1}^{n_i}\frac{(Y_{ij}-x^{\top}_{ij}\beta)^2}{v_i\sigma^2_{ij}(\tau)}+\sum_{j=1}^{n_i}\frac{s^2_{ij}(\alpha)}{\sigma^2_{ij}(\tau)}v_i \Bigg)\Bigg\}
\nonumber\\
&{}\qquad \times\frac{( \frac{\gamma}{\delta})^{\lambda}}{2K_\lambda(\delta\gamma)} v_i^{\lambda-1}\exp\left\{ -\frac{1}{2}\left( \frac{\delta^2}{v_i}+\gamma^2v_i \right)\right\}dv_i\Bigg]
\nonumber\\
&= \sum_{i=1}^N    \log\Bigg[ 
C_i(\al,\beta,\tau) \frac{( \gamma/\delta)^{\lambda}}{2K_\lambda(\delta\gamma)}
\int_0^{\infty} \exp\Bigg[-\frac{1}{2}\Bigg\{\frac{1}{v_i}\Bigg(\sum_{j=1}^{n_i}\frac{(Y_{ij}-x^{\top}_{ij}\beta)^2}{\sigma^2_{ij}(\tau)}+\delta^2\Bigg)
\nn\\
&{}\qquad 
+\Bigg(\sum_{j=1}^{n_i}\frac{s^2_{ij}(\alpha)}{\sigma^2_{ij}(\tau)}+\gamma^2\Bigg)v_i\Bigg\}\Bigg] v_i^{\lambda-1-\frac{n_i}{2}}dv_i\Bigg],
\nonumber
\end{align}
where
\begin{equation}
C_i(\al,\beta,\tau) := (2\pi)^{-\frac{n_i}{2}}\Bigg(\prod_{j=1}^{n_i} \sigma^2_{ij}(\tau)\Bigg)^{-\frac{1}{2}}\exp\Bigg(  \sum_{j=1}^{n_i}\frac{ s_{ij}(\alpha)(Y_{ij}-x^{\top}_{ij}\beta)}{\sigma^2_{ij}(\tau)}\Bigg).
\nonumber
\end{equation}
Making the change of variables $S_i^2 v_i/T_i = u_i$ with
\begin{align}
S_i=S_i(\al,\tau,\gam) &:= \sqrt{\gamma^2 + \sum_{j=1}^{n_i}\frac{s^2_{ij}(\alpha)}{\sigma^2_{ij}(\tau)}},
\nn\\
T_i=T_i(\beta,\al,\tau,\del,\gam) &:= 
S_i \,
\sqrt{\delta^2 + \sum_{j=1}^{n_i}\frac{1}{\sigma^2_{ij}(\tau)} (Y_{ij}-x^{\top}_{ij}\beta)^2},
\nonumber
\end{align}
we can continue as
\begin{align}
\ell_n(\theta)&= \sum_{i=1}^N    \log\left[ 
C_i(\al,\beta,\tau) \frac{( \frac{\gamma}{\delta})^{\lambda}}{2K_\lambda(\delta\gamma)}
\int_0^{\infty} \exp\left\{-\frac{1}{2}\left(\frac{S_i^2}{T_i u_i}\frac{T_i^2}{S_i^2}+\frac{T_i}{S_i^2}u_i S_i^2\right)\right\}
\frac{T_i^{\lambda-1-\frac{n_i}{2}}}{S_i^{2\left(\lambda-1-\frac{n_i}{2}\right)}}u_i^{\lambda-1-\frac{n_i}{2}}\frac{T_i}{S_i^2}du_i\right]
\nonumber\\
&= \sum_{i=1}^N \log\left[ 
C_i(\al,\beta,\tau) \frac{( \frac{\gamma}{\delta})^{\lambda}}{2K_\lambda(\delta\gamma)}
\int_0^{\infty} \exp\left\{ -\frac{T_i}{2}\left(\frac{1}{u_i}+u_i\right)\right\} 
\frac{T_i^{\lambda-\frac{n_i}{2}}}{S_i^{2\left(\lambda-\frac{n_i}{2}\right)}}u_i^{\lambda-1-\frac{n_i}{2}}du_{i}\right]
\nonumber\\
&= \sum_{i=1}^N \log\Biggl[ (2\pi)^{-\frac{n_i}{2}}\Bigg(\prod_{j=1}^{n_i} \sigma^2_{ij}(\tau)\Bigg)^{-1/2}\frac{(\gam/\del)^{\lambda}}{K_\lambda(\delta\gamma)}\exp\Bigg( \sum_{j=1}^{n_i}\frac{ s_{ij}(\alpha)}{\sigma^2_{ij}(\tau)}(Y_{ij}-x^{\top}_{ij}\beta) \Bigg)
\nonumber\\
&{}\qquad \times \frac{T_i^{\lambda-\frac{n_i}{2}}}{S_i^{2(\lambda-\frac{n_i}{2})}}\frac{1}{2}
\int_0^{\infty} \exp\left\{-\frac{T_i}{2}\left(\frac{1}{u_i}+u_i\right)\right\} u_i^{\lambda-1-\frac{n_i}{2}}du_{i}\Biggr]
\nonumber\\
&= \sum_{i=1}^N    \log\Biggl[ (2\pi)^{-\frac{n_i}{2}}\Bigg(\prod_{j=1}^{n_i} \sigma^2_{ij}(\tau)\Bigg)^{-1/2}\frac{(\gam/\del)^{\lambda}}{K_\lambda(\delta\gamma)}\exp\Bigg(  \sum_{j=1}^{n_i}\frac{ s_{ij}(\alpha)}{\sigma^2_{ij}(\tau)}(Y_{ij}-x^{\top}_{ij}\beta)\Bigg)
\left(\frac{T_i}{S_i^2}\right)^{\lambda-\frac{n_i}{2}}
K_{\lambda-\frac{n_i}{2}}(T_i)\Biggr].
\nonumber
\end{align}
This leads to the expression \eqref{log-LF}.

\subsection{Partial derivatives}
\label{sec_log-LF-partial.deri}

Recall the notation: $\ell_N(\theta)=\sumi \zeta_i(\theta)$, $A_i=A_i(\al,\tau,\gam)$ of \eqref{def_Ai}, and $B_i=B_i(\beta,\tau,\del)$ of \eqref{def_Bi}. Let
\begin{equation}
T_i' = T_i'(\theta) := R_{\lam-\frac{n_i}{2}}(A_i B_i) + \frac{\lam-\frac{n_i}{2}}{A_i B_i}
\nonumber
\end{equation}
for $R_\nu(t)$ defined by \eqref{def_R}. Then, we have the following expressions for the components of $\p_\theta\ell_N(\theta)$:
\rev{
\begin{align}
{\partial_\beta} \zeta_i(\theta)
&= -\sum_{j=1}^{n_i}\frac{s_{ij}(\alpha)}{\sigma^2_{ij}(\tau)}x_{ij} +\left(\lambda-\frac{n_i}{2}\right)\frac{1}{B_i}\partial_{\beta}B_i-T_i' A_i\partial_{\beta}B_i
\nonumber\\
&= -\sum_{j=1}^{n_i}\frac{s_{ij}(\alpha)}{\sigma^2_{ij}(\tau)}x_{ij} 
- \frac{1}{B_i} \left\{ \frac{1}{B_i}\left(\lambda-\frac{n_i}{2}\right) + T'_i A_i \right\}
\sum_{j=1}^{n_i}\frac{(Y_{ij}-x_{ij}^{\top}\beta)}{\sigma^2_{ij}(\tau)}x_{ij},
\nonumber\\
{\partial_\alpha} \zeta_i(\theta)
&= \sum_{j=1}^{n_i}\frac{\p_\al s_{ij}(\alpha)}{\sigma^2_{ij}(\tau)}(Y_{ij}-x_{ij}^{\top}\beta)-\left(\lambda-\frac{n_i}{2}\right)\frac{1}{A_i}\partial_\alpha A_i - T_i' B_i\partial_\alpha A_i
\nonumber\\
&=\sum_{j=1}^{n_i}\frac{\p_\al s_{ij}(\alpha)}{\sigma^2_{ij}(\tau)}(Y_{ij}-x_{ij}^{\top}\beta)
- \frac{1}{A_i} \left\{ \frac{1}{A_i}\left(\lambda-\frac{n_i}{2}\right) + T'_i B_i \right\}
\sum_{j=1}^{n_i}\frac{s_{ij}(\alpha)}{\sigma^2_{ij}(\tau)}\p_\al s_{ij}(\alpha),
\nonumber\\
{\partial_\tau}\zeta_i(\theta)
&= -\frac{1}{2}\sum_{j=1}^{n_i}\frac{\p_\tau (\sigma^2_{ij}(\tau))}{\sigma^2_{ij}(\tau)}
-\sum_{j=1}^{n_i}\frac{s_{ij}(\alpha)}{(\sigma^2_{ij}(\tau))^2}(Y_{ij}-x_{ij}^{\top}\beta)\p_\tau (\sigma^2_{ij}(\tau)) 
\nonumber\\
&{}\qquad +\left(\lambda-\frac{n_i}{2}\right)\left(\frac{\partial_\tau B_i}{B_i}-\frac{\partial_\tau A_i}{A_i}\right) - T_i'\partial_\tau (A_iB_i)
\nonumber\\
&= -\frac{1}{2}\sum_{j=1}^{n_i}\frac{\p_\tau (\sigma^2_{ij}(\tau))}{\sigma^2_{ij}(\tau)}
-\sum_{j=1}^{n_i}\frac{s_{ij}(\alpha)}{(\sigma^2_{ij}(\tau))^2}(Y_{ij}-x_{ij}^{\top}\beta)\p_\tau (\sigma^2_{ij}(\tau))
\nonumber\\
&{}\qquad -\frac12 \left(\lambda-\frac{n_i}{2}\right)
\Bigg(\frac{1}{B_i^2}\sum_{j=1}^{n_i}\frac{(Y_{ij}-x_{ij}^{\top}\beta)^2}{(\sigma^2_{ij}(\tau))^2}
\p_\tau (\sigma^2_{ij}(\tau)) -\frac{1}{A_i^2}\sum_{j=1}^{n_i}\frac{s^2_{ij}(\alpha)}{(\sigma^2_{ij}(\tau))^2}\p_\tau (\sigma^2_{ij}(\tau))\Bigg)
\nonumber\\
&{}\qquad + \frac12 T_i'
\Bigg(\frac{B_i}{A_i}\sum_{j=1}^{n_i}\frac{s^2_{ij}(\alpha)}{(\sigma^2_{ij}(\tau))^2}\p_\tau (\sigma^2_{ij}(\tau))
+\frac{A_i}{B_i}\sum_{j=1}^{n_i}\frac{(Y_{ij}-x_{ij}^{\top}\beta)^2}{(\sigma^2_{ij}(\tau))^2}\p_\tau (\sigma^2_{ij}(\tau))\Bigg),
\nonumber\\
{\partial_\lambda} \zeta_i(\theta)&= \log\left(\frac{\gamma}{\delta}\right) - \frac{ \partial_\lambda K_{\lambda}(\delta\gamma)}{ K_{\lambda}(\delta\gamma)}+\log B_i-\log A_i+\frac{ \partial_\lambda K_{\lambda-\frac{n_i}{2}}(A_iB_i)}{ K_{\lambda-\frac{n_i}{2}}(A_iB_i)},
\nonumber\\
{\partial_\delta} \zeta_i(\theta)
&= \gamma R_{\lambda}(\delta\gamma)+\left(\lambda-\frac{n_i}{2}\right)\frac{\delta}{B_i^2} - T_i' \frac{A_i}{B_i}\delta,
\nonumber\\
{\partial_\gamma} \zeta_i(\theta)
&= \frac{2\lambda}{\gamma}+\delta R_{\lambda}(\delta\gamma)-\left(\lambda-\frac{n_i}{2}\right)\frac{\gamma}{A_i^2} - T_i' \frac{B_i}{A_i}\gamma.
\nonumber
\end{align}
}



\medskip

\rev{
As for the second-order derivatives, for brevity, we write 
\begin{equation}
U_i = S_{\lambda-\frac{n_i}{2}}(A_iB_i)+\frac{1}{A_iB_i}R_{\lambda-\frac{n_i}{2}}(A_iB_i)-\frac{\lambda - \frac{n_i}{2}}{A^2_iB_i^2},
\nonumber
\end{equation}
for $R_\nu(t)$ and $S_\nu(t)$ defined by \eqref{def_R} and \eqref{def_S}. Further, let
\begin{equation}
L_\nu(z) := \frac{1}{K^2_{\nu}(z)}\left(\partial_{\lambda}K_{\nu-1}(z)K_{\nu}(z)-\partial_{\lambda}K_{\nu}(z)K_{\nu-1}(z)\right).
\nonumber
\end{equation}
Below we list the $21$ components of $\p_\theta^2\zeta_i(\theta)$, which were used to compute the confidence intervals and the one-step estimator;
the sizes of the matrices are not confusing, hence we are not taking care of them in notation and use the standard multilinear-form notation such as $(\p_\beta B_i) \otimes (\p_\al A_i) := \partial_\beta B_i\partial_\alpha^{\top} A_i \in \mbbr^{p_\beta}\otimes\mbbr^{p_\al}$.
\begin{align}
\p_\beta^2 \zeta_i(\theta)
&= \frac{\left(\lambda - \frac{n_i}{2}\right)}{B_i^2}\left(B_i \p_\beta^2 B_i -(\partial_\beta B_i)^{\otimes 2}\right) - T_i'
A_i \p_\beta^2 B_i -U_i A^2_i (\partial_\beta B_i)^{\otimes 2},
\nn\\
\partial_\beta \partial_\alpha \zeta_i(\theta)
&= -\sum_{j=1}^{n_i} \frac{1}{\sigma^2_{ij}(\tau)}(x_{ij} \otimes \p_\al s_{ij}(\alpha))
-(T_i' + U_i A_i B_i) \{(\p_\beta B_i) \otimes (\p_\al A_i)\},
\nn\\
\partial_\beta\partial_\tau \zeta_i(\theta)
&= \sum_{j=1}^{n_i}
\frac{s_{ij}(\alpha)}{\left(\sigma^2_{ij}(\tau)\right)^2}
\left(x_{ij} \otimes \p_\tau (\sigma^2_{ij}(\tau))\right)
+\frac{\left(\lambda-\frac{n_i}{2}\right)}{B_i^2}\left( B_i (\p_\tau\p_\beta B_i)
-(\partial_\beta B_i)\otimes (\partial_\tau B_i) \right) 
\nn\\
&{}\qquad - U_i A_i (\p_\beta B_i) \otimes (\p_\tau (A_i B_i)) - T_i' (\p_\beta B_i)\otimes (\p_\tau A_i) 
-T_i' A_i (\p_\beta \p_\tau B_i),
\nn\\
\partial_\beta\partial_\lambda\zeta_i(\theta)
&= A_i (\partial_\beta B_i) L_{\lambda-\frac{n_i}{2} }(A_iB_i),
\nn\\
\partial_\beta\partial_\delta\zeta_i(\theta)
&= -2\left(\lambda-\frac{n_i}{2}\right)\frac{\delta}{ B_i^{3}}\partial_\beta B_i
-U_iA_i^2\frac{\delta}{B_i}\partial_{\beta}B_i+T_i'A_i\frac{\delta}{B_i^2}\partial_\beta B_i,
\nn\\
\partial_\beta\partial_\gamma\zeta_i(\theta)
&= - \gamma U_i B_i \partial_{\beta}B_i - T_i'\frac{\gamma}{A_i}\partial_\beta B_i,
\nn\\
\partial_\alpha^2 \zeta_i(\theta)
&= \sum_{j=1}^{n_i}\frac{(y_{ij}-x_{ij}^{\top}\beta)}{\sigma^2_{ij}(\tau)} \p_\al^2 s_{ij}(\al)
-\frac{\left(\lambda-\frac{n_i}{2}\right)}{A_i^2}\left\{A_i \p_\al^2 A_i - (\partial_\alpha A_i)^{\otimes 2}\right\}
\nonumber\\
&{}\qquad -U_iB^2_i (\partial_{\alpha}A_i)^{\otimes 2} - T_i' B_i \partial_{\alpha}^2 A_i,
\nonumber\\
\partial_\alpha\partial_\tau\zeta_i(\theta)
&= -\sum_{j=1}^{n_i}\frac{(y_{ij}-x_{ij}^{\top}\beta)}{\sig_{ij}^4(\tau)}
\left\{(\p_\al s_{ij}(\al)) \otimes \p_\tau (\sigma^2_{ij}(\tau))\right\}
\nn\\
&{}\qquad - \frac{\left(\lambda-\frac{n_i}{2}\right)}{A_i^2}\left( A_i (\p_\al \p_\tau A_i) - (\p_\al A_i)\otimes (\p_\tau A_i)\right)
\nonumber\\
&{}\qquad -U_i B_i (\partial_\alpha A_i) \otimes (\partial_{\tau}(A_i B_i))
- T_i' \left( (\partial_\alpha A_i)\otimes (\partial_\tau B_i) + B_i \partial_\alpha\partial_\tau A_i)\right),
\nn\\
\partial_\alpha\partial_\lambda\zeta_i(\theta)
&= -2\frac{\partial_\alpha A_i}{A_i} - B_i (\partial_\alpha A_i) L_{\lambda-\frac{n_i}{2}}(A_iB_i),
\nn\\
\partial_\alpha\partial_\delta\ell_i(\theta)
&= -\del U_i A_i \partial_{\alpha}A_i - T_i'\frac{\delta}{B_i}\partial_\alpha A_i,
\nn\\
\partial_\alpha\partial_\gamma\zeta_i(\theta)
&= 2\left(\lambda-\frac{n_i}{2}\right)\frac{\gamma}{ A_i^{3}}\partial_\alpha A_i 
- U_iB^2_i\frac{\gamma}{A_i}\partial_{\alpha}A_i+T_i'B_i\frac{\gamma}{A_i^2}\partial_\alpha A_i,
\nn\\
\partial_\tau^2\zeta_i(\theta)
&= \sum_{j=1}^{n_i} \frac{s_{ij}(\alpha)(y_{ij}-x_{ij}^{\top}\beta)}{\sigma^8_{ij}(\tau)} 
\left(
2 \sigma^4_{ij(\tau)} \left(\p_\tau (\sigma^2_{ij}(\tau)) \right)^{\otimes 2}
- \sigma^4_{ij}(\tau) \p_\tau^2 (\sigma^2_{ij}(\tau))
\right)
\nonumber\\
&{}\qquad -\frac{1}{2} \sum_{j=1}^{n_i}\frac{\p_\tau^2 (\sigma^2_{ij}(\tau))\sigma^2_{ij}(\tau) - (\p_\tau (\sigma^2_{ij}(\tau)))^{\otimes 2}}{\sigma^4_{ij}(\tau)}
\nonumber\\
&{}\qquad +\left(\lambda-\frac{n_i}{2}\right)
\left[\frac{1}{B_i^2}\left(B_i (\partial_\tau^2 B_i) - (\partial_\tau B_i)^{\otimes 2}\right)
-\frac{1}{A_i^2}\left(A_i(\partial_\tau^2 A_i) - (\partial_\tau A_i)^{\otimes 2}\right)\right]
\nonumber\\
&{}\qquad - U_i\partial_{\tau}(A_iB_i) \otimes \partial_\tau(A_iB_i) - T_i' \partial_\tau^2(A_iB_i),
\nn\\
\partial_\tau\partial_\lambda\zeta_i(\theta)
&= -\frac{2\partial_\tau A_i}{A_i}-\partial_\tau (A_iB_i)L_{\lambda-\frac{n_i}{2}}(A_iB_i),
\nn\\
\partial_\tau\partial_\delta\zeta_i(\theta)
&= -2\left(\lambda-\frac{n_i}{2}\right)\frac{\delta}{ B_i^{3}}\partial_\tau B_i-U_iA_i\frac{\delta}{B_i}\partial_{\tau}(A_iB_i)+\frac{\delta T_i'}{B_i}\left(\frac{A_i}{B_i}\partial_\tau B_i-\partial_\tau A_i\right),
\nn\\
\partial_\tau\partial_\gamma\zeta_i(\theta)
&= 2\left(\lambda-\frac{n_i}{2}\right)\frac{\gamma}{ A_i^{3}}\partial_\tau A_i-U_iB_i\frac{\gamma}{A_i}\partial_{\tau}(A_iB_i)+\frac{\gamma T_i'}{A_i}\left(\frac{B_i}{A_i}\partial_\tau A_i-\partial_\tau B_i\right),
\nn\\
\partial^2_\lambda\zeta_i(\theta)
&= \frac{\partial^2_{\lambda}K_{\lambda-\frac{n_i}{2}}(A_iB_i)}{K_{\lambda-\frac{n_i}{2}}(A_iB_i)}-\left(\frac{\partial_{\lambda}K_{\lambda-\frac{n_i}{2}}(A_iB_i)}{K_{\lambda-\frac{n_i}{2}}(A_iB_i)}\right)^2-\frac{\partial^2_{\lambda}K_{\lambda-1}(\delta\gamma)}{K_{\lambda}(\delta\gamma)}+\left(\frac{\partial_{\lambda}K_{\lambda}(\delta\gamma)}{K_{\lambda}(\delta\gamma)}\right)^2,
\nn\\
\partial_\lambda\partial_\delta\zeta_i(\theta)
&= -A_i\frac{\delta}{B_i}L_{\lambda-\frac{n_i}{2}}(A_iB_i)+\gamma L_{\lambda}(\delta\gamma),
\nn\\
\partial_\lambda\partial_\gamma\zeta_i(\theta)
&= \frac{2}{\gamma} -\frac{2\gamma}{A_i^2}+\delta L_{\lambda}(\delta\gamma)-B_i\frac{\gamma }{A_i}L_{\lambda-\frac{n_i}{2}}(A_iB_i),
\nn\\
\partial_\delta^2\zeta_i(\theta)
&= \left(S_{\lambda}(\delta\gamma)+\frac{1}{\delta\gamma}R_{\lambda}(\delta\gamma)\right)\gamma^2+\left(\lambda-\frac{n_i}{2}\right)\left(\frac{1}{B_i^{2}}-\frac{2\delta^2}{B_i^{4}}\right)\nonumber\\
&{}\qquad -U_iA^2_i\frac{\delta^2}{B^2_i}-T_i'A_i\left(\frac{1}{B_i}-\frac{2\delta^2}{B_i^{3}}\right),
\nn\\
\partial_\delta\partial_\gamma\zeta_i(\theta)
&= 2R_{\lambda}(\delta\gamma)+\delta\gamma S_{\lambda}(\delta\gamma)-U_i\delta\gamma-T_i'\frac{\delta\gamma}{ A_iB_i},
\nn\\
\partial_\gamma^2\zeta_i(\theta)
&= -\frac{2\lambda}{\gamma^2}+\left(S_{\lambda}(\delta\gamma)+\frac{1}{\delta\gamma}R_{\lambda}(\delta\gamma)\right)\delta^2-\left(\lambda-\frac{n_i}{2}\right)\left(\frac{1}{A_i^{2}}-\frac{2\gamma^2}{A_i^{4}}\right)
\nn\\
&{}\qquad -U_iB^2_i\frac{\gamma^2}{A^2_i}-T_i'B_i\left(\frac{1}{A_i}-\frac{\gamma^2}{A_i^3}\right).
\nonumber
\end{align}
}


\end{document}